\documentclass[12pt]{amsart}
\usepackage{amsmath}
\usepackage{amssymb}
\usepackage{latexsym}

\newcommand{\Zint}{\mathbb {Z}}    
     
\newcommand{\halmos}{\rule{5pt}{5pt}}

\numberwithin{equation}{section}

\newtheorem{prop}{\bf Proposition}[section]
\newtheorem{thm}[prop]{\bf Theorem}

\newtheorem{cor}[prop]{\bf Corollary}

\begin{document}

\title[Heun equation {\rm IV}]
{The Heun equation and the Calogero-Moser-Sutherland system {\rm IV}: the Hermite-Krichever Ansatz}
\author{Kouichi Takemura}
\address{Department of Mathematical Sciences, Yokohama City University, 22-2 Seto, Kanazawa-ku, Yokohama 236-0027, Japan.}
\email{takemura@yokohama-cu.ac.jp}

\begin{abstract}
We develop a theory for the Hermite-Krichever Ansatz on the Heun equation.
As a byproduct, we find formulae which reduce hyperelliptic integrals to elliptic ones.
\end{abstract}

\maketitle

\section{Introduction}

Relations between ellptic integrals and hyperelliptic integrals have been studied from the 19th century.
For example Hermite \cite{Her} found the following formula
\begin{equation}
\int \frac{zdz}{\sqrt{(z^2-a)(8z^3-6az-b)}} = \frac{1}{2\sqrt{3}} \int \frac{dy}{\sqrt{y^3-3ay+b}}, \label{Hermitef}
\end{equation}
where $y=(2z^3-b)/(3(z^2-a))$.

In this paper we derive several formulae which reduce hyperelliptic integrals to elliptic ones. They are obtained by comparing two expressions of global monodromies of the Heun equation.
Here the Heun equation is the standard canonical form of a Fuchsian equation with four singularities, which is given by
\begin{equation}
\left( \! \left(\frac{d}{dw}\right) ^2 \! + \left( \frac{\gamma}{w}+\frac{\delta }{w-1}+\frac{\epsilon}{w-t}\right) \frac{d}{dw} +\frac{\alpha \beta w -q}{w(w-1)(w-t)} \right)\tilde{f}(w)=0
\label{Heun}
\end{equation}
with the condition 
\begin{equation}
\gamma +\delta +\epsilon =\alpha +\beta +1.
\label{Heuncond}
\end{equation}
It is also known that solving Heun equation corresponds to the spectral problem for a certain model of quantum mechanics which is called the $BC_1$ Inozemtsev model.
Set $\omega _0=0$ and $\omega_2=-\omega_1 -\omega _3$. The $BC_1$ Inozemtsev model is a one-particle quantum mechanics model whose Hamiltonian is given as
\begin{equation}
H= -\frac{d^2}{dx^2} + \sum_{i=0}^3 l_i(l_i+1)\wp (x+\omega_i),
\label{Ino}
\end{equation}
where $\wp (x)$ is the Weierstrass $\wp$-function with periods $(2\omega _1,2\omega _3 )$, $\omega _0, \omega _1, \omega _2, \omega_3$ are half-periods, and $l_i$ $(i=0,1,2,3)$ are coupling constants.
 This model is a one-particle version of the $BC_N$ Inozemtsev system \cite{Ino}, which is known to be the universal quantum integrable system with $B_N$ symmetry \cite{Ino,OOS}.
Let $f(x)$ be the eigenvector of the Hamiltonian $H$ whose eigenvalue is $E$, i.e.
\begin{equation}
(H-E) f(x)= \left( -\frac{d^2}{dx^2} + \sum_{i=0}^3 l_i(l_i+1)\wp (x+\omega_i)-E\right) f(x)=0.
\label{InoEF+}
\end{equation}
Then the coupling constants $l_0, l_1, l_2, l_3$ correpond to parameters $\alpha, \beta, \gamma, \delta, \epsilon$, a ratio of periods of elliptic function corresponds to a singular point $t$, an eigenvalue $E$ corresponds to an accessory $q$, and an eigenfunction $f(x)$ corresponds to a solution $\tilde{f}(w)$. For details see \cite{Smi,Tak1,Tak2,Tak3}.
For the case $l_0 \neq 0$ and $l_1=l_2=l_3=0$, Eq.(\ref{InoEF+}) is called the Lam\'e equation.

Hermite and Halphen investigated the Lam\'e equation by hypothesizing that solutions are expressed by an elliptic Baker-Akhiezer (Block) type function and Krichever \cite{Kri} described elliptic solutions to Kadomtset-Petviashvili equation by a Baker-Akhiezer type function. In this paper we investigate solutions to Eq.(\ref{InoEF+}) by considering the Hermite-Krichever Ansatz. In our situation, the Hermite-Krichever Ansatz asserts that the differential equation has solutions that are expressed as a finite series in the derivatives of an elliptic Baker-Akhiezer function, multiplied by an exponential function. More precisely, we are going to find solutions to Eq.(\ref{InoEF+}) of the form
\begin{align}
& f(x) = \exp \left( \kappa x \right) \left( \sum _{i=0}^3 \sum_{j=0}^{l_i-1} \tilde{b} ^{(i)}_j  \left( \frac{d}{dx} \right) ^{j} \Phi _i(x, \alpha ) \right) ,
\label{Lalpha+}
\end{align}
where $\Phi _i(x,\alpha )= \exp (\zeta( \alpha )x) \sigma (x+\omega _i -\alpha ) /\sigma (x+\omega _i )$ $(i=0,1,2,3)$. 

Treibich and Verdier \cite{TV} found and showed that, if $l_i \in \Zint _{\geq 0}$ for all $i \in \{ 0,1,2,3 \}$, then the potential $\sum_{i=0}^3 l_i(l_i+1)\wp (x+\omega_i)$ satisfies the stationary higher KdV equation and they constructed a theory of elliptic soliton following Krichever's idea \cite{Kri}; while Gesztesy, Weikard \cite{GW,Wei}, Smirnov \cite{Smi} and the author \cite{Tak1,Tak3} obtained further results on this subject.
Thus, the function $ \sum_{i=0}^3 l_i(l_i+1)\wp (x+\omega_i)$ is called the Treibich-Verdier potential; and is closely related with a hyperelliptic curve $\nu^2 = -Q(E)$ where $Q(E)$ is a polynomial in $E$ which is determined for each $l_0, l_1, l_2, l_3 \in \Zint $. For example, if the eigenvalue $E$ satisfies $Q(E)=0$, then the equation has a doubly periodic eigenfunction up to signs which corresponds to the Heun polynomial. (See \cite{GW,Smi,Tak1,Tak2}.) In \cite{Tak3}, global monodomies of Eq.(\ref{InoEF+}) with the condition $l_0, l_1, l_2, l_3 \in \Zint $ were calculated and they are expressed by hyperelliptic integral (see Proposition \ref{thm:conj3}).

Belokolos, Eilbeck, Enolskii, Kostov and Smirnov studied the covering map from the hyperelliptic curve $\nu^2 = -Q(E)$ to the ellitic curve $\wp '(\alpha)^2=4 \wp (\alpha)^3-g_2 \wp (\alpha) -g_3$ and they obtained relations among variables $E, \alpha $ and $ \kappa$ in Eq.(\ref{Lalpha+}) for the case $l_0=1,2,3,4,5$, $l_1=l_2=l_3=0$, the case $l_0=2, l_1=1, l_2=0, l_3=0$ and the case $l_0=2, l_1 =1 ,l_2 =1, l_3=0$ (see \cite{BE} and the reference therein). By considering the covering map they found transformation formulae like Eq.(\ref{Hermitef}) that reduce hyperelliptic integrals to elliptic ones case by case. 
On the other hand, Maier found a pattern of the covering maps for the case of Lam\'e equation (i.e. the case $l_0 \neq 0$, $l_1=l_2=l_3=0$) and conjectured formulae \cite[Conjecture L]{Mai} by introducing the notions ``twisted Lam\'e polynomials'' and ``theta-twisted Lam\'e polynomials'' .

In this paper we justify and develop the Hermite-Krichever Ansatz on the Heun equation without an advanced algebraic geometry technique for the case $l_0, l_1, l_2, l_3 \in \Zint _{\geq 0}$. Note that results on the Bethe Ansatz and monodromy formulae in terms of hyperelliptic integrals obtained in \cite{Tak1,Tak3} play important roles in our approach. 
As a result, the monodromies of Eq.(\ref{InoEF+}) are expressed by elliptic integrals.
To study the Heun equation by the Hermite-Krichever Ansatz, we need to consider the covering map in detail.
For this purpose, we introduce twisted Heun and theta-twisted Heun polynomials that are based on Maier's ideas, and we obtain theorems that support the Maier's conjectures.
By comparing two expressions of monodromies of Eq.(\ref{InoEF+}),
we establish transformation formulae between elliptic integrals of the first and second kinds, and hyperelliptic integrals of the first and second kinds. Hence, the mysteries of the elliptic-hyperelliptic integral formulae are unveiled by the monodromies.
For the case $l_0=2$, $l_1=l_2=l_3=0$, we obtain Eq.(\ref{Hermitef}) as a transformation formula between elliptic integrals of first kind and hyperelliptic integrals of first kind. The formula for the second kind is written as
\begin{equation}
\int \frac{(2z^2-a)dz}{\sqrt{(z^2-a)(8z^3-6az-b)}} -\frac{1}{3}\sqrt{\frac{8z^3-6az-b}{z^2-a}} = \frac{1}{2\sqrt{3}} \int \frac{ydy}{\sqrt{y^3-3ay+b}}, \label{Hermites}
\end{equation}
where $y=(2z^3-b)/(3(z^2-a))$.

The Hermite-Krichever Ansatz would be applicable to the spectral problem of the $BC_1$ Inozemtsev model, because the monodromy is expressed in terms of an elliptic integral by applying the Hermite-Krichever Ansatz, and it is closely related with the boundary condition of the model.

This paper is organized as follows. In section \ref{sec:HKA}, we justify the Hermite-Krichever Ansatz on the Heun equation by applying results on the Bethe Ansatz and integral representation of solutions.
In section \ref{sec:hee}, we obtain hyperelliptic-elliptic reduction formulae by comparing two expressions of monodromies.
In section \ref{sec:covering}, we investigate zeros and poles of the covering map.
In section \ref{sec:tHttH}, we introduce twisted Heun and theta-twisted Heun polynomials, and obtain formulae which support Maier's conjectures.
In section \ref{sec:exa}, we give examples which cover the cases that the genus of the related hyperelliptic integral is less than or equal to three. Thus we obtain hyperelliptic-elliptic reduction formulae explicitly for more than 20 cases.

Throughtout this paper we assume $l_0, l_1, l_2, l_3 \in \Zint_{\geq 0}$ and $(l_0, l_1 ,l_2 ,l_3) \neq (0,0,0,0)$.

\section{Hermite-Krichever Ansatz} \label{sec:HKA}

In this section we review results on Bethe Ansatz and integral representation of solutions which are obtained in \cite{Tak1}, and apply them to justify the Hermite-Krichever Ansatz.

Fix the eigenvalue $E$ of the Hamiltonian $H$ (see Eq.(\ref{Ino})) and consider the second-order differential equation
\begin{equation}
(H-E) f(x)= \left( -\frac{d^2}{dx^2} + \sum_{i=0}^3 l_i(l_i+1)\wp (x+\omega_i)-E\right) f(x)=0.
\label{InoEF}
\end{equation}

Let $h(x)$ be  the product of any pair of the solutions to Eq.(\ref{InoEF}). Then the function $h(x)$ satisfies the following third-order differential equation:
\begin{align}
& \left( \frac{d^3}{dx^3}-4\left( \sum_{i=0}^3 l_i(l_i+1)\wp (x+\omega_i)-E\right)\frac{d}{dx}-2\left(\sum_{i=0}^3 l_i(l_i+1)\wp '(x+\omega_i)\right) \right) h (x)=0.
\label{prodDE} 
\end{align}

It is known that Eq.(\ref{prodDE}) has a nonzero doubly periodic solution for all $E$ if $l_i \in \Zint _{\geq 0}$ $(i=0,1,2,3)$.
\begin{prop} \cite[Proposition 3.5]{Tak1} \label{prop:prod}
If $l_0, l_1, l_2, l_3 \in \Zint _{\geq 0}$, then equation (\ref{prodDE}) has a nonzero doubly periodic solution $\Xi (x,E)$, which has the expansion
\begin{equation}
\Xi (x,E)=c_0(E)+\sum_{i=0}^3 \sum_{j=0}^{l_i-1} b^{(i)}_j (E)\wp (x+\omega_i)^{l_i-j},
\label{Fx}
\end{equation}
where the coefficients $c_0(E)$ and $b^{(i)}_j(E)$ are polynomials in $E$, they do not have common divisors and the polynomial $c_0(E)$ is monic.
We set $g=\deg_E c_0(E)$. Then the coefficients satisfy $\deg _E b^{(i)}_j(E)<g$ for all $i$ and $j$.
\end{prop}

We can derive the integral formula for the solution $\Lambda(x,E)$ to Eq.(\ref{InoEF}) in terms of the doubly periodic function $\Xi(x,E)$, which is obtained in \cite{Tak1}. Set
\begin{align}
 & Q(E)=  \Xi (x,E)^2\left( E- \sum_{i=0}^3 l_i(l_i+1)\wp (x+\omega_i)\right) +\frac{1}{2}\Xi (x,E)\frac{d^2\Xi (x,E)}{dx^2}-\frac{1}{4}\left(\frac{d\Xi (x,E)}{dx} \right)^2. \label{const}
\end{align}
It is shown in \cite{Tak1} that $Q(E)$ is independent of $x$. Thus $Q(E)$ is a monic polynomial in $E$ of degree $2g+1$, which follows from the expression for $\Xi (x,E)$ given by Eq.(\ref{Fx}). The following proposition on the integral representation of solutions is obtained in \cite{Tak1}:
\begin{prop} \cite[Proposition 3.7]{Tak1}
Let $\Xi (x,E)$ be the doubly periodic function defined in Proposition \ref{prop:prod} and $Q(E)$ be the monic polynomial defined in Eq.(\ref{const}).
Then the function 
\begin{equation}
\Lambda (x,E)=\sqrt{\Xi (x,E)}\exp \int \frac{\sqrt{-Q(E)}dx}{\Xi (x,E)}
\label{eqn:Lam}
\end{equation}
is a solution to the differential equation (\ref{InoEF}).
\end{prop}

We consider the case $Q(E)=0$. 
It follows from Eq.(\ref{eqn:Lam}) that, $\Lambda (x,E) ^2= \Xi(x,E)$. 
By considering zeros and poles, we obtain
\begin{equation}
\Lambda (x,E) ^2= \Xi(x,E) =\frac{ C\prod_{j=1}^{( \beta _0 +\beta _1 +\beta _2 +\beta _3 )/2}(\wp(x)-\wp(t_j))^2}{(\wp(x)-e_1)^{\beta _1}(\wp(x)-e_2)^{\beta_2}(\wp(x)-e_3)^{\beta _3}} \label{Fxtj0}
\end{equation}
where $\beta _i \in \{ l_i, -l_i -1 \}$ $(i=0,1,2,3)$, $C, t_1,\dots $ are constants. From Eq.(\ref{eq:Leg}), the function 
\begin{equation}
\tilde{\Lambda }(x,E) = \frac{\prod_{j=1}^{( \beta _0 +\beta _1 +\beta _2 +\beta _3 )/2} \sigma(x-t_j)\sigma(x+t_j)}{\sigma(x)^{\beta _0}\sigma_1(x)^{\beta _1}\sigma_2(x)^{\beta _2}\sigma_3(x)^{\beta _3}}, 
\label{eq:tilL0}
\end{equation}
is a solution to Eq.(\ref{InoEF}). We set $\alpha =  \sum_{k=1}^3 \beta _k\omega_k $ for the case $Q(E) = 0$.
Then we have 
\begin{equation} 
\tilde{\Lambda }(x+2\omega _k,E)   = 
\left\{ \begin{array}{ll}
\tilde{\Lambda }(x,E) & (\alpha \equiv 0 \mbox{ (mod }2\omega_1 \Zint \oplus 2\omega_3 \Zint  )) \\
\exp (-2\eta _k \alpha +2\omega _k \zeta (\alpha ) ) \tilde{\Lambda }(x,E) & (\alpha \not\equiv 0 \mbox{ (mod }2\omega_1 \Zint \oplus 2\omega_3 \Zint  ))
\end{array} \right.
\label{ellint0} 
\end{equation}
for $k=1,2,3$.

Now we consider the case $Q(E) \neq 0$. Then it is known that the functions $\Lambda (x,E)$ and $\Lambda (-x,E)$ are linearly independent.

Set $l=l_0+l_1+l_2 +l_3$. The function $\Xi (x,E)$ is an even doubly periodic function which may have poles of degree $2l_i$ at $x=\omega_i$ $(i=0,1,2,3)$.
Therefore we have the following expression;
\begin{equation}
\Xi (x,E)=\frac{b_0^{(0)} (E)\prod_{j=1}^{l}(\wp(x)-\wp(t_j))}{(\wp(x)-e_1)^{l_1}(\wp(x)-e_2)^{l_2}(\wp(x)-e_3)^{l_3}} \label{Fxtj}
\end{equation}
for some values $t_1,\dots ,t_l$. It is shown in \cite{Tak1} that, if $Q(E) \neq 0$, then the values $t_j$, $-t_j$, $\omega_i$ $($mod $2\omega_1 \Zint \oplus 2\omega_3 \Zint )$ $(j=1,\dots ,l,\; i=0,1,2,3)$ are mutually distinct. Set $z=\wp(x)$ and $z_j=\wp(t_j)$. From Eq.(\ref{const}) we have $\left( \frac{d\Xi (x,E)}{dz}\right)^2 _{z=z_j}=\frac{-4Q(E)}{\wp'(t_j)^2}$. We fix the signs of $t_j$ by taking 
\begin{equation}
\left( \frac{d\Xi (x,E)}{dz}\right) _{z=z_j}=-\frac{2\sqrt{-Q(E)}}{\wp'(t_j)}.
\label{signpptj}
\end{equation}
Let $\zeta (x)$ be the Weierstrass zeta function, $\sigma(x)$ be the Weierstrass sigma function and $\sigma_i(x)$ $(i=1,2,3)$ be the Weierstrass co-sigma functions. (See appendix.) 
If we put $2\sqrt{-Q(E)}/\Xi (x,E)$ into partial fractions, it is seen that
\begin{align}
& \frac{2\sqrt{-Q(E)}}{\Xi (x,E)}=\sum_{j=1}^l \frac{\frac{2\sqrt{-Q(E)}}{\left( \frac{d\Xi}{dz}\right) _{z=z_j}}}{\wp (x)-\wp (t_j)}  = \sum_{j=1}^l \left(\zeta(x-t_j)-\zeta(x+t_j)+2\zeta(t_j) \right) \nonumber ,
\end{align}
It follows that
\begin{align}
& \Lambda (x,E)= \label{BA} \\
& =\sqrt{\frac{b_0^{(0)} (E)\prod_{j=1}^{l}(\wp(x)-\wp(t_j))}{\prod_{k=1}^3 (\wp(x)-e_k)^{l_k}}} \exp\left( \frac{1}{2}\sum_{j=1}^l \left( \log \sigma(t_j+x)-\log \sigma (t_j-x) -2x\zeta(t_j)\right) \right) \nonumber \\
& = \frac{\sqrt{b_0^{(0)} (E)}\prod_{j=1}^l \sigma(x+t_j)}{\sigma(x)^{l_0}\sigma_1(x)^{l_1}\sigma_2(x)^{l_2}\sigma_3(x)^{l_3}\prod_{j=1}^l \sigma(t_j)}\exp \left(-x\sum_{i=1}^l \zeta(t_j)\right). \nonumber 
\end{align}
For the case $Q(E)\neq 0$, we set
\begin{align}
& \tilde{\Lambda }(x,E) = \frac{\prod_{j=1}^l \sigma(x+t_j)}{\sigma(x)^{l_0}\sigma_1(x)^{l_1}\sigma_2(x)^{l_2}\sigma_3(x)^{l_3}\prod_{j=1}^l \sigma(t_j)}\exp \left(-x\sum_{i=1}^l \zeta(t_j)\right), 
\label{eq:tilL}
\end{align}
Then the formula
\begin{equation}
\Lambda (x,E) = \sqrt{b_0^{(0)} (E)} \tilde{\Lambda }(x,E)
\label{LamLamtilde}
\end{equation}
is derived.

We establish the validity of the Hermite-Krichever Ansatz by using values $t_1, \dots , t_l$.
More precisely we show that, if $l_0 , l_1,l_2,l_3 \in \Zint _{\geq 0}$, then an eigenfunction of the Hamiltionian $H$ with every eigenvalue $E$ is expressed as in the form of Theorem \ref{thm:alpha}. 
We set $\alpha =-\sum_{j=1}^{l} t_j+ \sum_{k=1}^3 l_k\omega_k$ and $\kappa = \zeta (\sum_{j=1}^{l} t_j - \sum_{k=1}^3 l_k\omega_k) - \sum_{j=1}^{l} \zeta (t_j)+ \sum_{k=1}^3 l_k \eta_k$, where $\eta _k =\zeta (\omega _k)$ $(k=1,2,3)$. Since the set $\{ -t_j \}_{j=1,\dots ,l }$ is the set of zeros of $\tilde{\Lambda }(x,E)$ and the position of the zeros and the poles are doubly-periodic, $\pm \alpha $ is determined uniquely mod $2\omega_1 \Zint \oplus 2\omega_3 \Zint $ and $\kappa $ is determined uniquely.
From Eqs.(\ref{periods}, \ref{eq:sigmai}) it follows that
\begin{align} 
& \tilde{\Lambda }(x+2\omega _k,E)  \label{ellint} \\
& = \exp \left( 2\eta _k \left( \sum_{j=1}^l t_j - \sum_{k'=1}^3 l_{k'}\omega_{k'} \right)  -2 \omega _k \left( \sum _{j=1}^l \zeta (t_j ) -\sum_{k'=1}^3 l_{k'}\eta _{k'}  \right)  \right)  \tilde{\Lambda }(x,E) \nonumber \\
&  = \exp (-2\eta _k \alpha +2\omega _k \zeta (\alpha ) +2 \kappa \omega _k ) \tilde{\Lambda }(x,E) \nonumber
\end{align}
for $k=1,2,3$.
We set
\begin{equation}
\Phi _i(x,\alpha )= \frac{\sigma (x+\omega _i -\alpha ) }{ \sigma (x+\omega _i )} \exp (\zeta( \alpha )x), \quad \quad (i=0,1,2,3).
\label{Phii}
\end{equation}
Then we have 
\begin{equation}
\left( \frac{d}{dx} \right) ^{j} \Phi _i(x+2\omega _k , \alpha ) = \exp (-2\eta _k \alpha +2\omega _k\zeta (\alpha )) \left( \frac{d}{dx} \right) ^{j} \Phi _i(x, \alpha )
\label{ddxPhiperiod}
\end{equation}
for $i=0,1,2,3$, $k=1,2,3$ and $j \in \Zint _{\geq 0}$.

The following theorem asserts that eigenfunctions of the Hamitonian $H$ are expressed in the form of the Hermite-Krichever Ansatz.
\begin{thm} \label{thm:alpha}
Assume $l_0, l_1, l_2 , l_3 \in \Zint_{\geq 0}$ and set $l= l_0+l_1+l_2+l_3$. Let $\{ -t_j \} _{j=1, \dots ,l}$ be the set of zeros of $\tilde{\Lambda }(x,E)$ appeared in Eq.(\ref{eq:tilL}). If $\alpha = -\sum_{j=1}^{l} t_j+ \sum_{k=1}^3 l_k\omega_k \not\equiv 0$ $($mod $2\omega_1 \Zint \oplus 2\omega_3 \Zint)$, then we have
\begin{align}
& \tilde{\Lambda} (x,E) = \exp \left( \kappa x \right) \left( \sum _{i=0}^3 \sum_{j=0}^{l_i-1} \tilde{b} ^{(i)}_j \left( \frac{d}{dx} \right) ^{j} \Phi _i(x, \alpha ) \right)
\label{Lalpha}
\end{align}
for some values $\tilde{b} ^{(i)}_j$ $(i=0,\dots ,3, \: j= 0,\dots ,l_i-1)$. The values $\alpha $ and $\kappa $ are expressed as
\begin{equation}
 \wp (\alpha ) =\frac{P_1 (E)}{P_2 (E)}, \; \; \; \wp ' (\alpha ) =\frac{P_3 (E)}{P_4 (E)} \sqrt{-Q(E)} , \; \; \kappa  =\frac{P_5 (E)}{P_6 (E)} \sqrt{-Q(E)},
\label{P1P6}
\end{equation}
where $P_1 (E) ,\dots ,P_6(E)$ are polynomials in $E$.

If $\alpha = -\sum_{j=1}^{l} t_j+ \sum_{k=1}^3 l_k\omega_k \equiv 0$ $($mod $2\omega_1 \Zint \oplus 2\omega_3 \Zint)$, then we have
\begin{align}
& \tilde{\Lambda} (x,E) = \exp \left( \bar{\kappa } x \right) \left( \bar{c} +\sum _{i=0}^3 \sum_{j=0}^{l_i-2} \bar{b} ^{(i)}_j \left( \frac{d}{dx} \right) ^{j} \wp (x+\omega _i) +\sum_{k=1}^3 \bar{c}_k \frac{\wp '(x)}{\wp (x)-e_k} \right)
\label{Lalpha0}
\end{align}
for some values $\bar{c}$, $\bar{c}_k$ $(k=1,2,3)$ and $\bar{b} ^{(i)}_j$ $(i=0,\dots ,3, \: j= 0,\dots ,l_i-2)$.
\end{thm}
\begin{proof}
Assume $\alpha \not \equiv 0$. From Eq.(\ref{ellint}) and Eq.(\ref{ddxPhiperiod}), the functions $\exp(\kappa x) \left( \frac{d}{dx} \right) ^{j} \Phi _i(x , \alpha ) $ ($i=0,1,2,3$, $j \in \Zint _{\geq 0}$)
and the function $\tilde{\Lambda }(x,E)$ have the same periodicities.
By subtracting the functions $\exp(\kappa x) \left( \frac{d}{dx} \right) ^{j} \Phi _i(x , \alpha ) $ from the function $\tilde{\Lambda }(x,E)$ to erase poles, we obtain an holomorphic function that has the same periods as $\Phi _0(x , \alpha ) $ and it must be zero.
Hence we obtain the expression (\ref{Lalpha}). For the case $\alpha \equiv 0$, we change $\{ t_j \} _{j=1}^l $ by $t_1 \rightarrow t_1 +\alpha $ and $t_j \rightarrow t_j$ $(j=2,\dots ,l)$. Then we can set $\alpha = -\sum_{j=1}^{l} t_j+ \sum_{k=1}^3 l_k\omega_k =0$. By setting $\bar{\kappa } =- \sum_{j=1}^{l} \zeta (t_j)+ \sum_{k=1}^3 l_k \eta_k $ we obtain the expression (\ref{Lalpha0}).

Next we investigate the values $\wp (\alpha )$, $\wp ' (\alpha )$ and $\kappa $ for the case $\alpha \not \equiv 0$.
The functions $\wp (\sum_{j=1}^l t_j - \sum_{k=1}^3 l_{k}\omega_{k} )$, $\wp '(\sum_{j=1}^l t_j - \sum_{k=1}^3 l_{k}\omega_{k})$ and $(\kappa =) \zeta (\sum_{j=1}^{l} t_j - \sum_{k=1}^3 l_{k}\omega_{k}) - \sum_{j=1}^{l} \zeta (t_j) + \sum_{k=1}^3 l_{k}\eta_{k}$ are doubly-periodic in variables $t_1 ,\dots ,t_l$.
Hence by applying addition formulae of elliptic functions and considering the parity of functions $\wp (x)$, $\wp '(x)$ and $\zeta (x)$, we obtain the expression
\begin{align}
& \wp (\sum_{j=1}^l t_j - \sum_{k=1}^3 l_{k}\omega_{k}) = \sum_{j_1<j_2<\dots <j_m\atop{m:\mbox{\scriptsize{ even}}}}f^{(1)}_{j_1, \dots ,j_m}(\wp(t_1), \dots ,\wp (t_l)) \wp'(t_{j_1}) \dots \wp'(t_{j_l}), \\
& \wp ' (\sum_{j=1}^l t_j - \sum_{k=1}^3 l_{k}\omega_{k}) = \sum_{j_1<j_2<\dots <j_m\atop{m:\mbox{\scriptsize{ odd}}}}f^{(2)}_{j_1, \dots ,j_m}(\wp(t_1), \dots ,\wp (t_l)) \wp'(t_{j_1}) \dots \wp'(t_{j_l}), \nonumber \\
& \zeta (\sum_{j=1}^{l} t_j - \sum_{k=1}^3 l_{k}\omega_{k}) - \sum_{j=1}^{l} \zeta (t_j) + \sum_{k=1}^3 l_{k}\eta_{k} \nonumber \\
& = \sum_{j_1<j_2<\dots <j_m\atop{m:\mbox{\scriptsize{ odd}}}}f^{(3)}_{j_1, \dots ,j_m}(\wp(t_1), \dots ,\wp (t_l)) \wp'(t_{j_1}) \dots \wp'(t_{j_l}), \nonumber
\end{align}
where $f^{(k)}_{j_1, \dots ,j_m}(x_1, \dots , x_l)$ $(k=1,2,3)$ are rational functions in $x_1 ,\dots x_l$.
By applying Eq.(\ref{signpptj}), the function  $\wp '(t_j)/\sqrt{-Q(E)}$ is expressed as a rational function in $E$ and $\wp (t_j)$.
Hence $\wp (\sum_{j=1}^l t_j - \sum_{k=1}^3 l_{k}\omega_{k})$, $\wp '(\sum_{j=1}^l t_j - \sum_{k=1}^3 l_{k}\omega_{k})/\sqrt{-Q(E)}$ and $(\zeta (\sum_{j=1}^{l} t_j- \sum_{k=1}^3 l_{k}\omega_{k}) - \sum_{j=1}^{l} \zeta (t_j) + \sum_{k=1}^3 l_{k}\eta_{k})/\sqrt{-Q(E)}$ are expressed as rational functions in the variable $\wp(t_1), \dots ,\wp (t_l)$ and $E$.
They are symmetric in $\wp(t_1), \dots ,\wp (t_l)$.

From Eq.(\ref{Fxtj}) we have
\begin{equation}
b_0^{(0)} (E)\prod_{j=1}^{l}(\wp(x)-\wp(t_j))= \Xi (x,E) (\wp(x)-e_1)^{l_1}(\wp(x)-e_2)^{l_2}(\wp(x)-e_3)^{l_3}.
\end{equation}
Hence the elementary symmetric functions $ \sum _{j_1<\dots <j_{l'}} \wp(t_{j_1})\dots \wp (t_{j_{l'}})$ ($l'=1,\dots ,l$) are expressed as rational functions in $E$.
Therefore $\wp (\sum_{j=1}^l t_j - \sum_{k=1}^3 l_{k}\omega_{k})$, $\wp '(\sum_{j=1}^l t_j - \sum_{k=1}^3 l_{k}\omega_{k})/\sqrt{-Q(E)}$ and $(\zeta (\sum_{j=1}^{l} t_j- \sum_{k=1}^3 l_{k}\omega_{k}) - \sum_{j=1}^{l} \zeta (t_j) + \sum_{k=1}^3 l_{k}\eta_{k})/\sqrt{-Q(E)}$ are expressed as rational functions in $E$.
\end{proof}

\begin{prop} \label{alphastar}
We have $\alpha \rightarrow 0$ $($mod $2\omega_1 \Zint \oplus 2\omega_3 \Zint )$ as $E \rightarrow \infty$.
\end{prop}
\begin{proof}
We define the rational function $\tilde{\Xi }(z,E)$ to satisfy 
\begin{equation}
\tilde{\Xi }(\wp (x),E) =\Xi (x,E)/E^g.
\end{equation}
From Proposition \ref{prop:prod}, we have
\begin{equation}
\tilde{\Xi }(z,E) = 1 +\frac{1}{E} \frac{P(z,1/E)}{(z-e_1)^{l_1}(z-e_2)^{l_2}(z-e_3)^{l_3}}
\end{equation}
where $P(x,y)$ is a polynomial in $x$ and $y$.
Let $k \in \{ 1,2,3 \}$. If $E$ is sufficiently large and $|z-e_k| = |E|^{-\frac{1}{2(l_k+2)}}$, then we have $|\tilde{\Xi} (z,E) -1| <|E|^{-1/2}$ and $|\frac{\partial }{\partial z}\tilde{\Xi} (z,E) |<|E|^{-1/2}$. Hence we have
\begin{equation}
\frac{1}{2\pi \sqrt{-1}} \oint _{|z-e_k| = |E|^{-\frac{1}{2(l_k+2)}}} \frac{\frac{\partial }{\partial z}\tilde{\Xi} (z,E)}{\tilde{\Xi} (z,E)} \rightarrow 0
\label{intzeros}
\end{equation}
as $|E| \rightarrow \infty$.
Since the l.h.s. of Eq.(\ref{intzeros}) is equal to the number of zeros inside $|z-e_k| < |E|^{-\frac{1}{2(l_k+2)}}$ minus the number of poles inside $|z-e_k| < |E|^{-\frac{1}{2(l_k+2)}}$, if $|E|$ is sufficiently large, then the function $\tilde{\Xi} (z,E)$ has $l_k$ zeros in variable $z$ inside the circle $|z-e_k| = |E|^{-\frac{1}{2(l_k+2)}}$.
It follows from Proposition \ref{prop:prod} that the degree of $P(z,1/E)$ in $z$ is $l_0+l_1+l_2+l_3$. Hence it is shown similarly that, if $|E|$ is sufficiently large, then the function $\tilde{\Xi} (z,E)$ has $l_0$ zeros in variable $z$ outside the circle $|z| = |E|^{\frac{1}{2(l_0+2)}}$.
Therefore, if $|E| \rightarrow \infty$, then $l_k$ zeros of the function $\tilde{\Xi} (z,E)$ in $z$ tend to $e_k$ $(k=1,2,3)$ and $l_0$ zeros tend to infinity.

Since the set $\{ \wp (-t_j) \} _{j=1,\dots , l_0+l_1+l_2+l_3}$ coincide with the set of zeros of $\tilde{\Xi} (z,E)$ in $z$, $l_i$ elements of $\{-t_j \} _{j=1,\dots , l_0+l_1+l_2+l_3}$ tends to $\omega _i$ $(i=0,1,2,3)$.
Hence we have $\alpha = \sum _{j=1}^{l_0+l_1+l_2+l_3} (-t_j) + \sum_{k=1}^3 l_k\omega _k \rightarrow 0$ (mod $2\omega_1 \Zint \oplus 2\omega_3 \Zint$) as $E\rightarrow \infty$.
Thus we obtain the proposition.
\end{proof}

\section{Hyperelliptic-elliptic reduction formulae} \label{sec:hee}

We obtain hyperelliptic-elliptic reduction formulae by comparing two expressions of monodromies.
One is Eq.(\ref{ellint}), and the other one is the monodromy formula in terms of hyperelliptic integral that was obtained in \cite{Tak3}.
\begin{prop} \cite[Theorem 3.7]{Tak3} \label{thm:conj3} 
Assume $l_i \in \Zint_{\geq 0}$ ($i=0,1,2,3$). Let  $k \in \{ 1,2,3\}$, $q_k \in \{0,1\}$ and $E_0$ be the eigenvalue such that $\Lambda (x+2\omega _k,E_0)=(-1)^{q_k} \Lambda (x,E_0)$. Then
\begin{equation} 
\Lambda (x+2\omega _k,E)=(-1)^{q_k} \Lambda (x,E) \exp \left( -\frac{1}{2} \int_{E_0}^{E}\frac{ \int_{0+\varepsilon }^{2\omega _k+\varepsilon }\Xi (x,\tilde{E})dx}{\sqrt{-Q(\tilde{E})}} d\tilde{E}\right)
\label{analcontP}
\end{equation}
with $\varepsilon $ a constant so determined as to avoid passing through the poles in the integration. 
\end{prop}

Since the function $\wp (x+\omega_i)^n $ is written as a linear combination of the functions $\left( \frac{d}{dx} \right) ^{2j} \wp (x+\omega_i)$ $(j=0, \dots ,n)$, the function $\Xi (x,E)$ can be expressed as
\begin{equation}
\Xi (x,E)=c(E)+\sum_{i=0}^3 \sum_{j=0}^{l_i-1 } a^{(i)}_j (E)\left( \frac{d}{dx} \right) ^{2j} \wp (x+\omega_i),
\label{FFx}
\end{equation}
Set
\begin{equation}
a(E)=\sum _{i=0}^3 a^{(i)} _0 (E).
\label{polaE}
\end{equation}
From Proposition \ref{thm:conj3} we have
\begin{equation}
\tilde{\Lambda }(x+2\omega _k,E)=(-1)^{q_k} \tilde{\Lambda }(x,E) \exp \left( -\frac{1}{2} \int_{E_0}^{E}\frac{ -2\eta _k a(\tilde{E}) +2\omega _k c(\tilde{E}) }{\sqrt{-Q(\tilde{E})}} d\tilde{E}\right) .
\label{hypellint}
\end{equation}

Set 
\begin{align}
& A=-2\alpha -\int_{E_0}^{E}\frac{a(\tilde{E})}{\sqrt{-Q(\tilde{E})}} d\tilde{E} , \label{eq:defA} \\
& B=2\zeta (\alpha )+2\kappa + \int_{E_0}^{E}\frac{c(\tilde{E}) }{\sqrt{-Q(\tilde{E})}} d\tilde{E} .
\end{align}
By comparing Eq.(\ref{ellint}) and Eq.(\ref{hypellint}) for the cases $k=1$ and $k=3$, we have 
\begin{align}
& \eta _1 A +\omega _1 B=\pi \sqrt{-1} (q_1+2n_1),\\
& \eta _3 A +\omega _3 B=\pi \sqrt{-1} (q_3+2n_3),
\end{align}
for some integers $n_1$ and $n_3$.
By the Legendre's relation $\eta _1 \omega _3 - \eta _3 \omega _1 =\pi\sqrt{-1}/2$
(see Eq.(\ref{eq:Leg})), we have
\begin{align}
& A/2 =(q_1+2n_1)\omega_3 -(q_3+2n_3)\omega _1 , \label{eq:Ao3o1}\\
& -B/2 =(q_1+2n_1)\eta_3 -(q_3+2n_3)\eta _1.
\end{align}
Hence
\begin{align}
& \alpha +\frac{1}{2}\int_{E_0}^{E}\frac{a(\tilde{E})}{\sqrt{-Q(\tilde{E})}} d\tilde{E} =-(q_1+2n_1)\omega_3 +(q_3+2n_3)\omega _1, \label{alpE} \\
& \zeta (\alpha )+\kappa + \frac{1}{2}\int_{E_0}^{E}\frac{c(\tilde{E}) }{\sqrt{-Q(\tilde{E})}} d\tilde{E} =-(q_1+2n_1)\eta_3 +(q_3+2n_3)\eta _1. \label{alpEz}
\end{align}
We set $\xi =\wp (\alpha )$. By Proposition \ref{alphastar} and the relation $\int (1/\wp'(\alpha )) d\xi =\int d\alpha$, we have
\begin{equation}
\int _{\infty} ^{\xi } \frac{d \tilde{\xi }}{\sqrt{4 \tilde{\xi } ^3-g_2 \tilde{\xi } -g_3}} = \alpha = -\frac{1}{2} \int _{\infty}^{E} \frac{a(\tilde{E})}{\sqrt{-Q(\tilde{E})}}d\tilde{E}
\label{alpint}
\end{equation}
Note that $Q(E)$ is a polynomial with degree $2g+1$, and $a(E)$ is a polynomial with degree $g$. Hence Eq.(\ref{alpint}) represents a formula which reduces a hyperelliptic integral of the first kind to an elliptic integral of the first kind. The transformation of variables is given by $\xi =P_1(E) /P_2 (E)$ for polynomials $P_1 (E)$ and $P_2 (E)$ (see Eq.(\ref{P1P6})). In sections \ref{sec:covering} and \ref{sec:tHttH} we investigate the transformation of variables in detail. 
Let $\alpha _0$ be the value $\alpha $ evaluated at $E=E_0$. It follows from Eqs.(\ref{eq:defA}, \ref{eq:Ao3o1}) that $\alpha _0 =-(q_1+2n_1)\omega_3 +(q_3+2n_3)\omega _1$ and 
\begin{align}
& \alpha - \alpha _0 +\frac{1}{2}\int_{E_0}^{E}\frac{a(\tilde{E})}{\sqrt{-Q(\tilde{E})}} d\tilde{E} =0.
\label{aaaQ}
\end{align}
If $\alpha _0 \equiv 0$ (mod $2\omega _1 \Zint \oplus 2\omega _3 \Zint$), then $\zeta (\alpha -\alpha _0 )= \zeta (\alpha ) + (q_1+2n_1)\eta_3  - (q_3+2n_3)\eta _1$, which follows from $\zeta ( \alpha +2\omega _k) =\zeta (\alpha )+ 2\eta _k$ $(k=1,3)$. Combining with Eqs.(\ref{alpEz}, \ref{aaaQ}) we have
\begin{align}
& \kappa = -\frac{1}{2} \int  _{E_0}^{E} \frac{c(\tilde{E})}{\sqrt{-Q(\tilde{E})}}d\tilde{E} + \zeta \left(  \frac{1}{2}\int_{E_0}^{E}\frac{a(\tilde{E})}{\sqrt{-Q(\tilde{E})}} d\tilde{E} \right) . \label{kap0}
\end{align}
If $\alpha _0 \not\equiv 0$ (mod $2\omega _1 \Zint \oplus 2\omega _3 \Zint$), then $\zeta (\alpha _0 )= -(q_1+2n_1)\eta_3 +(q_3+2n_3)\eta _1$ and 
\begin{align}
 \kappa & = -\frac{1}{2} \int _{E_0}^{E} \frac{c(\tilde{E})}{\sqrt{-Q(\tilde{E})}}d\tilde{E} - \zeta (\alpha) + \zeta (\alpha _0) \label{kap123} \\
 & = -\frac{1}{2} \int  _{E_0}^{E} \frac{c(\tilde{E})}{\sqrt{-Q(\tilde{E})}}d\tilde{E} + \int _{\wp (\alpha _0)} ^{\xi } \frac{ \tilde{\xi } d \tilde{\xi }}{\sqrt{4\tilde{\xi }^3-g_2 \tilde{\xi }-g_3}} . \nonumber
\end{align}
Note that $Q(E)$ is a polynomial with degree $2g+1$, and $c(E)$ is a polynomial with degree $g+1$. Hence Eq.(\ref{kap123}) represents a formula which reduces a hyperelliptic integral of the second kind to an elliptic integral of the second kind.
The following proposition presents asymptotic behaviour of $\wp (\alpha )$ and $\kappa $ as $E \rightarrow \infty$.
\begin{prop} \label{prop:asym}
As $E \rightarrow \infty$, we have $\alpha \sim \frac{1}{2\sqrt{- E}} \sum_{i=0}^3 l_i(l_i+1)$, $\wp (\alpha ) \sim -4E/(\sum_{i=0}^3 l_i(l_i+1))^2$ and $\kappa \sim \sqrt{- E} (1-2/(\sum_{i=0}^3 l_i(l_i+1)))$.
\end{prop}
\begin{proof}
Set $\Xi (x,E)= \sum_{i=0}^g f_{i}(x) E^i$. By definition we have  $f_g(x)=1$. Substituting into (\ref{prodDE}) and comparing coefficients of $E^g$, we have $f'_{g-1} (x)=v'(x)/2$, where $v(x)= \sum_{i=0}^3 l_i(l_i+1) \wp (x+\omega _i)$. Hence $f_{g-1}(x)= v(x)/2 +C$ for some constants $C$.
Therefore we have $c(E) \sim E^g $, $a(E) \sim E^{g-1}\sum_{i=0}^3 l_i(l_i+1)/2 $  as $E\rightarrow \infty$. From Eq.(\ref{const}) we have $Q(E) \sim E^{2g+1}$ as $E\rightarrow \infty$. Thus
\begin{align}
& \int _{\infty}^{E} \frac{a(\tilde{E})}{\sqrt{-Q(\tilde{E})}}d\tilde{E} \sim \int _{\infty}^{E} \frac{\tilde{E}^{g-1}\sum_{i=0}^3 l_i(l_i+1)/2}{\sqrt{-1}\tilde{E}^{(2g+1)/2}}d\tilde{E} \sim -\sum_{i=0}^3 l_i(l_i+1) /\sqrt{-E} \label{inta}\\
&  \int _{E_0}^{E} \frac{c(\tilde{E})}{\sqrt{-Q(\tilde{E})}}d\tilde{E} \sim \int _{E_0}^{E} \frac{\tilde{E}^{g}}{\sqrt{-1}\tilde{E}^{(2g+1)/2}}d\tilde{E} \sim -2\sqrt{-E} \label{intc}
\end{align}
as $E\rightarrow \infty$. By combining Eq.(\ref{inta}) with Eq.(\ref{alpint}) we have 
$\alpha \sim \frac{1}{2\sqrt{-E}} \sum_{i=0}^3 l_i(l_i+1)$ as $E \rightarrow \infty $.
It is known that $\wp (x) \rightarrow x^{-2}$ and $\zeta (x) \rightarrow x^{-1}$ as $x \rightarrow 0$. Hence we have $\wp (\alpha )\sim -4E/(\sum_{i=0}^3 l_i(l_i+1))^2$ and $\zeta (\alpha )\sim 2\sqrt{-E}/(\sum_{i=0}^3 l_i(l_i+1))$ as $E\rightarrow \infty$.
By combining with Eqs.(\ref{kap0}, \ref{kap123}) and Eq.(\ref{intc}) we have $\kappa \sim \sqrt{- E} (1-2/(\sum_{i=0}^3 l_i(l_i+1)))$.
\end{proof}

\section{Covering map} \label{sec:covering}

Let $T$ be the elliptic curve whose basis periods are $(2\omega _1, 2\omega_3)$ and $V_{Q(E)}$ be the hyperelliptic curve defined by $\nu ^2=-Q(E)$. 
In Theorem \ref{thm:alpha}, a map $\pi : \: V_{Q(E)} \rightarrow T$ which is called covering or covering map in \cite{BE,Mai} is defined.
Let $\alpha \in T$ be the image of $(E,\nu )\in V_{Q(E)}$, i.e. $\alpha =\pi (E ,\nu )$.
From Theorem \ref{thm:alpha} and Proposition \ref{prop:asym} we obtain the following expression of the covering map:
\begin{equation}
\wp (\alpha )- e_k =-\frac{4N_k(E)}{(\sum_{i=0}^3 l_i(l_i+1))^2D_k(E)}, \quad (k=1,2,3),
\label{wpek}
\end{equation}
where $N_k(E)$ and $D_k(E)$ are monic polynomials such that $\deg _E N_k(E)= 1 +\deg _E D_k(E)$.
In this section we investigate zeros of the polynomials $N_k(E)$ and $D_k(E)$ $(k=1,2,3)$. 

The polynomials $N_k(E)$ and $D_k(E)$ $(k=1,2,3)$ can be calculated as the proof of Theorem \ref{thm:alpha}, but practically it is hard to do. In \cite{BE}, the value $\wp (\alpha )$ is calculated by expanding functions as 
\begin{align}
& \wp (x)= \frac{1}{x^2} +\frac{g_2 x^2}{20}+ \dots , \\
& \frac{\Phi_0 (x,\alpha )}{\sigma (-\alpha )}= \frac{1}{x} - \frac{1}{2}\wp(\alpha )x + \frac{1}{6} \wp '(\alpha ) x^2+ \dots ,
\end{align}
substituting them into Eqs.(\ref{InoEF}, \ref{Lalpha}) and comparing coefficients of $x^j$ $(j=-2,-1, \dots )$, althought it is not still effective.
Maier found a pattern of the covering maps for the case of Lam\'e equation and conjectured formulae of the covering maps \cite[Conjecture L]{Mai} by introducing the notions ``twisted Lam\'e polynomials'' and ``theta-twisted Lam\'e polynomials''.
In this paper we calculate the polynomials $N_k(E)$ and $D_k(E)$ $(k=1,2,3)$ by developing Maier's ideas.

Fix the value $E$. Then the function $\Lambda (x,E)$ and its zeros $\{-t_j\} _{j=1,\dots ,l}$ are determined, although we have ambiguity of choosing  $\Lambda (-x,E)$ and its zeros $\{t_j\} _{j=1,\dots ,l}$. Hence the value $\alpha (= -\sum_{j=1}^{l} t_j+ \sum_{k=1}^3 l_k\omega_k )$ (mod $2\omega_1 \Zint \oplus 2\omega_3 \Zint$) is determined up to the sign $\pm$. Note that this ambiguity corresponds to choosing $\nu $ or $-\nu $ on the hyperelliptic curve $\nu ^2=-Q(E)$. The value $\wp (\alpha )$ is determined uniquely.

Let $p$ be an element of $\{ 0,1,2,3 \}$ such that $\sum _{i=0}^3 l_i \omega _i \equiv \omega _p$ (mod $2\omega_1 \Zint \oplus 2\omega_3 \Zint$).
Assume that $\alpha \equiv \omega _p$ (mod $2\omega_1 \Zint \oplus 2\omega_3 \Zint$). Let $\wp _k (x)$ $(k=1,2,3)$ be co-$\wp$ function defined in Eq.(\ref{eq:sigmai}). It follows from Theorem \ref{thm:alpha} and Eqs.(\ref{periods}, \ref{rel:sigmai}) that the function $\tilde{\Lambda } (x,E) \exp (-\bar{\kappa } x) \wp _1 (x)^{l_1} \wp _2 (x)^{l_2} \wp _3 (x)^{l_3} $ is doubly-periodic. Hence the eigenfunction $\tilde{\Lambda } (x,E)$ is written as 
\begin{equation}
\tilde{\Lambda} (x,E) = \exp \left( \bar{\kappa } x \right) 
\frac{\left( \sum_{j=0}^{\hat{l}^{(0)}} \bar{a}_j (z-e_2) ^j \right) + \sqrt{(z-e_1)(z-e_2)(z-e_3)}\left( \sum_{j=0}^{\check{l}^{(0)}} \bar{b}_j (z-e_2)^j \right) }{(z-e_1)^{l_1/2}(z-e_2)^{l_2/2}(z-e_3)^{l_3/2}},
\label{Lomega0}
\end{equation}
where $z=\wp (x)$ and $\hat{l}^{(0)}$ (resp. $\check{l}^{(0)}$) is the greatest integer that is less than or equal to $(l_0+l_1+l_2+l_3)/2$ (resp. $(l_0+l_1+l_2+l_3-3)/2$). 
Let $A_p$ (resp. $B_p$) be the set of eigenvalues $E$ of the operator $H$ such that the eigenfunction $\tilde{\Lambda} (x,E) $ is written as the form (\ref{Lomega0}) and satisfies $\bar{\kappa }=0$ (resp. $\bar{\kappa }\neq 0$).

Next assume $\alpha \equiv \omega _{p'}$ (mod $2\omega_1 \Zint \oplus 2\omega_3 \Zint$) for $p' \in \{0,1,2,3\} \setminus \{ p\} $. Let $k$ be an element of $\{1,2,3\}$ such that $\omega _{p'} \equiv \omega _{p}+ \omega _{k}$ (mod $2\omega_1 \Zint \oplus 2\omega_3 \Zint$). 
It follows from Theorem \ref{thm:alpha} and Eqs.(\ref{periods}, \ref{rel:sigmai}) that the function $\tilde{\Lambda } (x,E) \exp (-\kappa  x) \wp _1 (x)^{l_1} \wp _2 (x)^{l_2} \wp _3 (x)^{l_3} $ has the same periodicity as the functions $\wp _k (x)$ and $\wp _{k'} (x)\wp _{k''} (x)$, where $k', k''$ are determined as $\{k,k',k''\}=\{1,2,3\}$. Hence the eigenfunction $\tilde{\Lambda } (x,E)$ is written as 
\begin{equation}
 \tilde{\Lambda} (x,E) = \exp \left( \kappa x \right) 
\frac{\sqrt{z-e_k} \left( \sum_{j=0}^{\hat{l}} a_j (z-e_k)^j \right) + \sqrt{(z-e_{k'})(z-e_{k''})}\left( \sum_{j=0}^{\check{l}} b_j (z-e_k)^j \right) }{(z-e_1)^{l_1/2}(z-e_2)^{l_2/2}(z-e_3)^{l_3/2}},
\label{Lomega}
\end{equation}
where $z=\wp (x)$ and $\hat{l}$ (resp. $\check{l}$) is the greatest integer that is less than or equal to $(l_0+l_1+l_2+l_3-1)/2$ (resp. $(l_0+l_1+l_2+l_3-2)/2$). 
Let $A_{p'}$  (resp. $B_{p'}$) be the set of eigenvalues $E$ of the operator $H$ such that the eigenfunction $\tilde{\Lambda} (x,E) $ is written as the form (\ref{Lomega}) and satisfies $\kappa =0$  (resp. $\kappa \neq 0$).

For the sets $A_i$, $B_i$ $(i=0,1,2,3)$, we have 
\begin{prop} \label{prop:AiBi}
(i) The eight sets $A_i$, $B_i$ $(i=0,1,2,3)$ are pairwise distinct.\\
(ii) The set $\cup _{i=0}^3 A_i$ coincides with the set of zeros of the polynomial $Q(E)$.\\
(iii) Assume that $l_0+l_1+l_2+l_3$ is even (resp. odd) and $E \in B_p$. Then we have $\bar{a}_{\hat{l}^{(0)}} \neq 0$ (resp. $\bar{b}_{\check{l}^{(0)}} \neq 0$).\\
(iv) Assume that $l_0+l_1+l_2+l_3$ is even (resp. odd) and $E \in (\cup _{i=0}^3 B_i) \setminus B_p$. Then we have $b_{\check{l}} \neq 0$ (resp. $a_{\hat{l}} \neq 0$).
\end{prop}
\begin{proof}
It follows from Theorem \ref{thm:alpha} that the values $\kappa^2$ and $\wp (\alpha)$ are determined uniquely for each $E$. If $E \in A_i \cup B_i$, then we have $\alpha \equiv \omega _i$ (mod $2\omega_1 \Zint \oplus 2\omega_3 \Zint$) and $\wp (\alpha )= \left\{ \begin{array}{ll} \infty & (i=0) ,\\ e_i & (i=1,2,3). \end{array} \right. $ Thus the sets $A_i$, $B_i$ $(i=0,1,2,3)$ are separated by the values of $\kappa^2$ and $\wp (\alpha)$. Hence (i) is proved.

Let ${\mathcal F}$ be the space spanned by meromorphic doubly periodic functions up to signs, namely
\begin{align}
& {\mathcal F}=\bigoplus _{\epsilon _1 , \epsilon _3 =\pm 1 } {\mathcal F} _{\epsilon _1 , \epsilon _3 }, \label{spaceF} \\
& {\mathcal F} _{\epsilon _1 , \epsilon _3 }=\{ f(x) \mbox{: meromorphic }| f(x+2\omega_1)= \epsilon _1 f(x), \; f(x+2\omega_3)= \epsilon _3 f(x) \} .
\end{align}
Then it is shown in \cite[Theorem 3.8]{Tak1} that Eq.(\ref{InoEF}) has a non-zero solution in the space ${\mathcal F}$ if and only if the value $E'$ satisfies the equation $Q(E')=0$.
Suppose $E' \in A_i$ for some $i$. Then the value $\bar{\kappa }$ or $\kappa $ is equal to zero. It follows from Eq.(\ref{Lomega0}) or Eq.(\ref{Lomega}), that there exists a non-zero eigenfunction $\tilde{\Lambda }(x,E')$ satisfying $\tilde{\Lambda }(x,E') \in {\mathcal F} _{\epsilon _1 , \epsilon _3 }$ for some $\epsilon _1 , \epsilon _3 \in \{ -1,1\}$. Hence we have $Q(E') =0$ and $\cup _{i=0}^3 A_i \subset \{E\: |\: Q(E)=0\}$.

Conversely suppose that $Q(E')=0$. Then we can take the eigenvalue $E_0$ in Proposition \ref{thm:conj3} as $E_0=E'$. From Eq.(\ref{alpE}) we have $\alpha \in \omega_1 \Zint \oplus \omega _3 \Zint$. Hence the function $\tilde{\Lambda }(x,E') $ is expressed as Eq.(\ref{Lomega0}) or Eq.(\ref{Lomega}). If $\bar{\kappa }\neq 0$ or $\kappa \neq 0$, then the functions $\tilde{\Lambda }(x,E')$ and $\tilde{\Lambda }(-x,E')$ are a basis of solutions to Eq.(\ref{InoEF}) and they are not doubly-periodic up to signs. It contradicts the existence of a solution in the space ${\mathcal F}$. Thus we have $\bar{\kappa }= 0$ or $\kappa =0$. Therefore we have $\cup _{i=0}^3 A_i \supset \{E\: |\: Q(E)=0\}$ and obtain (ii).

We show (iii). Assume that $\bar{a}_{\hat{l}^{(0)}} = 0$ (resp. $\bar{b}_{\check{l}^{(0)}} = 0$). Then the degree of the pole at $x=0$ of the functions $\tilde{\Lambda} (x,E)$ and $\tilde{\Lambda} (-x,E)$ are less than $l_0$. Since $\bar{\kappa } \neq 0$, the functions $\tilde{\Lambda} (x,E)$ and $\tilde{\Lambda} (-x,E)$ are linearly independent. Hence the degree of the pole at $x=0$ of all solutions to Eq.(\ref{InoEF}) must be less than $l_0$. But it contradicts that the exponents of solutions to Eq.(\ref{InoEF}) are $-l_0$ and $l_0 +1$. Hence we obtain (iii). (iv) is proved similarly.
\end{proof}

It follows from definitions of $A_i$, $B_i$ $(i=0,1,2,3)$ that, if $E \in \cup _{i=0}^3 A_i \cup B_i$, then a solution to Eq.(\ref{InoEF}) is written as Eq.(\ref{Lomega0}) or Eq.(\ref{Lomega}).
Conversely, it is shown that, if a solution to Eq.(\ref{InoEF}) is written as Eq.(\ref{Lomega0}) or Eq.(\ref{Lomega}), then we have $E \in \cup _{i=0}^3 A_i \cup B_i$. More precisely, 
\begin{prop} \label{prop:Lalp}
(i) If a solution to Eq.(\ref{InoEF}) is written as Eq.(\ref{Lomega0}), then we have $E \in A_p \cup B_p$.\\
(ii) Let $k \in \{1,2,3\}$ and $p'$ be the element of $\{ 0,1,2,3 \}$ such that $\omega _{p'} \equiv \omega _{p}+\omega _k$ (mod $2\omega_1 \Zint \oplus 2\omega_3 \Zint$). If a solution to Eq.(\ref{InoEF}) is written as Eq.(\ref{Lomega}), then we have $E \in A_{p'} \cup B_{p'}$.
\end{prop}
\begin{proof}
We show (ii).
Let $\tilde{\Lambda } ^a (x,E)$ be a solution to Eq.(\ref{InoEF}) written as Eq.(\ref{Lomega}). Let $S _E$ be the set of all solutions to Eq.(\ref{InoEF}) for a fixed eigenvalue $E$. Then $S_E$ is a two dimensional vector space.
Let $M_{k'}$ $(k'=1,3)$ be the transformations obtained by the analytic continuation $x \rightarrow x+2\omega _{k'}$. Then the actions of $M_{k'}$ preserve the space $S_E$ and it is shown in \cite[\S 3.2]{Tak1} that the transformations $M_{k'}$ do not depend on the choice of paths and $M_1 M_3= M_3 M_1$.
It follows from Eqs.(\ref{eq:sigmai}, \ref{rel:sigmai}) that 
\begin{align} 
& M_{k'} (\tilde{\Lambda }^a (x,E) ) = \exp (-2\eta _{k'} \omega _{p'} +2\omega _{k'} \tilde{\kappa })\tilde{\Lambda }^a (x,E) \label{eq:Mka}
\end{align}
for some value $\tilde{\kappa }$ and $k'=1,3$.

Assume that $Q(E) \neq 0$. Then it is shown in \cite{Tak1} the functions $\tilde{\Lambda} (x,E)$ and $\tilde{\Lambda} (-x,E)$ defined in Eq.(\ref{eq:tilL}) form a basis of the space $S _E$.
Since the function $\tilde{\Lambda} (x,E)$ is written as Eq.(\ref{ellint}), we have\begin{align} 
& M_{k'} (\tilde{\Lambda }(x,E)) = \exp (-2\eta _{k'} \alpha +2\omega _{k'} \chi )\tilde{\Lambda }(x,E) \label{eq:Mk'}
\end{align}
for some value $\chi$ and $k'=1,3$. It follows from changing variable $x \rightarrow -x-2\omega _{k'}$ that 
\begin{align} 
& M_{k'} (\tilde{\Lambda }(-x,E)) = \exp (2\eta _{k'} \alpha -2\omega _{k'} \chi )\tilde{\Lambda }(-x,E) .\label{eq:Mk'-}
\end{align}
From the assumption $Q(E) \neq 0$ and \cite[Theorem 3.8]{Tak1} we have $S _E \cap \tilde{F} = \{ 0\}$. Hence we obtain that $M_1$ or $M_3$ have eigenvalues of neither $1$ nor $-1$. Let $M_{\tilde{k}}$ be the one that has an eigenvalue of neither $1$ nor $-1$. Then the eigenvalues of $M_{\tilde{k}}$ are distinct. On the other hand the function $\tilde{\Lambda }^a (x,E) $
is also an eigenvector of the operators $M_{\tilde{k}}$ (see Eq.(\ref{eq:Mka})). Hence the function $\tilde{\Lambda }^a (x,E) $ is written as $C\tilde{\Lambda} (x,E)$ or $C\tilde{\Lambda} (-x,E)$ for some non-zero constant $C$. 
From Eqs.(\ref{eq:Mka}, \ref{eq:Mk'}) or Eqs.(\ref{eq:Mka}, \ref{eq:Mk'-}) we have
\begin{align}
& \pm (-2\eta _1 \alpha +2\omega _1 \chi )= -2\eta _1 \omega _{p'} +2\omega _1 \tilde{\kappa } +2\pi \sqrt{-1}n_1 \label{eq:n1} \\
& \pm (-2\eta _3 \alpha +2\omega _3 \chi )= -2\eta _3 \omega _{p'} +2\omega _3 \tilde{\kappa } +2\pi \sqrt{-1}n_3, \label{eq:n3}
\end{align}
for integers $n_1$, $n_3$. It follows that
\begin{equation}
 (\pm \alpha -\omega _{p'}) (-2\eta_1 \omega _3 + 2\eta_3 \omega _1 ) = 2\pi \sqrt {-1} (n_1\omega _3 -n_3 \omega _1 ) .
\end{equation}
From the Legendre's relation $\eta_1 \omega _3 - \eta_3 \omega _1 = \pi \sqrt {-1}/2$ (see Eq.(\ref{eq:Leg})) and the relation $\omega _{p'} \equiv \pm \omega _{p'}$ (mod $2\omega_1 \Zint \oplus 2\omega_3 \Zint$) we have $\alpha \equiv \omega _{p'}$ (mod $2\omega_1 \Zint \oplus 2\omega_3 \Zint$). Therefore $E \in B_{p'}$.

Assume that $Q(E) =0$. Then the set of the eigenvalues of the operators $M_{k'}$ $(k'=1,3)$ on the space $S_E$ is either $\{1 \}$ or $\{ -1 \}$, and the function $\tilde{\Lambda} (x,E)$ in Eq.(\ref{eq:tilL0}) is an eigenvector of the operators $M_{k'}$ $(k'=1,3)$ with an eigenvalue given in Eq.(\ref{ellint0}). Hence the eigenvalue of the operator $M_{k'}$ admits the expression as in Eq.(\ref{eq:Mk'}).
On the other hand the function $\tilde{\Lambda }^a (x,E)$ is also an eigenvector of the operators $M_{k'}$ $(k' =1,3)$ with an eigenvalue given in Eq.(\ref{eq:Mka}).
By comparing eigenvalues we obtain Eqs.(\ref{eq:n1}, \ref{eq:n3}) and that $\alpha \equiv \omega _{p'}$ (mod $2\omega_1 \Zint \oplus 2\omega_3 \Zint$). Therefore $E \in A_{p'}$.

(i) is proved similarly.
\end{proof}

The value $\wp (\alpha )- e_k$ $(k=1,2,3)$ is described by using the sets $A_i$, $B_i$ $(i=0,k)$.
\begin{prop} \label{thm:paek}
Let $m(E' )$ be the integer such that $a(E) =(E-E' )^{m(E' )} \tilde{a}(E)$ and $ \tilde{a}(E' ) \neq 0$, and $\tilde{m}(E' )$ be the integer such that $Q(E) =(E-E' )^{\tilde{m}(E' )} \tilde{Q}(E)$ and $ \tilde{Q}(E' ) \neq 0$ .
We set
\begin{align}
& D(E)= \prod _{E' \in A_0} (E-E' )^{2+2m(E' )-\tilde{m}(E' ) } \prod _{E' \in B_0} (E-E' )^{2+2m(E' )} ,
\label{eq:DE}\\
& N_k(E)= \prod _{E' \in A_{k}} (E-E' )^{2+2m(E' )-\tilde{m}(E' )} \prod _{E' \in B_{k}} (E-E' )^{2+2m(E' )}, \quad (k=1,2,3). \label{eq:NkE}
\end{align}
Then we have
\begin{equation}
\wp (\alpha )= e_k -\frac{4N_k(E)}{(\sum_{i=0}^3 l_i(l_i+1))^2D(E)}
\label{eq:wpaND}
\end{equation}
and $\deg _E N_k(E)= 1 +\deg _E D(E)$ for $k=1,2,3$.
\end{prop}
\begin{proof}
It follows from Theorem \ref{thm:alpha} and Proposition \ref{prop:asym} that $\wp (\alpha ) -e_k$ $(k=1,2,3)$ is expressed as Eq.(\ref{wpek}) such that $D_k(E)$ and $N_k(E)$ are monic polynomials, coprime and $\deg _E N_k(E)= 1 +\deg _E D_k(E)$.
We show that polynomials $N_k(E)$ and $D_k(E)$ are expressed as Eqs.(\ref{eq:DE}, \ref{eq:NkE}).

We consider poles and zeros of the rational function $\wp(\alpha )-e_k$ in variable $E$. Assume that $\alpha \not\equiv 0,\omega _k$ (mod $2\omega _1 \Zint \oplus 2\omega _3 \Zint$). Then $\wp (\alpha )-e_k$ is non-zero and non-infinity. Hence the eigenvalue $E'$ such that $\alpha \not\equiv 0, \omega _k$ does not satisfy $D(E')=0$ nor $N_k(E')=0$.

Next assume that $E'$ is a zero of the rational function $\wp(\alpha )-e_k$ in variable $E$. We write $\wp(\alpha )-e_k = (E-E')^m R (E)$ ($m \in \Zint _{\geq 1}$, $R (E') \neq 0, \infty$). Then $E'$ is an eigenvalue of the Hamiltonian that satisfies  $\alpha \equiv \omega _k$ (mod $2\omega _1 \Zint \oplus 2\omega _3 \Zint$). From Eq.(\ref{alpint}) it follows that 
\begin{equation}
\int _{e_k} ^{\wp (\alpha )} \frac{d\xi}{\sqrt{4\xi^3-g_2\xi-g_3}} = -\frac{1}{2} \int _{E'}^{E} \frac{a(\tilde{E})}{\sqrt{-Q(\tilde{E})}}d\tilde{E} .
\label{alpinte2}
\end{equation}
As $E \rightarrow E'$, we have
\begin{align}
& \int _{e_k} ^{\wp (\alpha )} \frac{d\xi}{\sqrt{4\xi^3-g_2\xi-g_3}} \sim \frac{\sqrt{\wp(\alpha )-e_k}}{\sqrt{(e_k-e_{k'})(e_k-e_{k''})}}=C_1(E-E')^{m/2} ,\\
& -\frac{1}{2} \int _{E'}^{E} \frac{a(\tilde{E})}{\sqrt{-Q(\tilde{E})}}d\tilde{E}
=  \int _{E'}^{E} (\tilde{E} -E')^{m(E')-\tilde{m}(E')/2}\tilde{R}(\tilde{E})d\tilde{E} \\
& \sim C_2 (E -E')^{1+ m(E')-\tilde{m}(E')/2}, \nonumber 
\end{align}
for some functions $\tilde{R}(E)$ and non-zero constants $C_1, C_2$.
Hence we have $m=2+2m(E')-\tilde{m}(E')$. 
It follows from Proposition \ref{prop:AiBi} that, if $E' \in B_k$, then $\tilde{m}(E')=0$. Therefore, if $E' \in A_k$ (resp. $E' \in B_k$), then the multiplicity of zeros of the polynomial $N_k(E)$ at $E=E'$ is $2+2m(E')-\tilde{m}(E')$ (resp. $2+2m(E')$).

We finally assume that $E'$ is a pole of the rational function $\wp(\alpha )-e_k$ in variable $E$. Then we have $\alpha \rightarrow 0$ (mod $2\omega _1 \Zint \oplus 2\omega _3 \Zint$) as $E\rightarrow E'$. We write $\wp(\alpha )-e_k = R (E)/(E-E')^m $ ($m \in \Zint _{\geq 1}$, $R (E') \neq 0, \infty$). From Eq.(\ref{alpint}) it follows that 
\begin{equation}
\int _{\infty} ^{\wp (\alpha )} \frac{d\xi}{\sqrt{4\xi^3-g_2\xi-g_3}} = -\frac{1}{2} \int _{E'}^{E} \frac{a(\tilde{E})}{\sqrt{-Q(\tilde{E})}}d\tilde{E} .
\label{alpintinfty}
\end{equation}
As $E \rightarrow E'$, we have
\begin{align}
& \int _{\infty } ^{\wp (\alpha )} \frac{d\xi}{\sqrt{4\xi^3-g_2\xi-g_3}} \sim \frac{-1}{\sqrt{\wp(\alpha )}}=\bar{C}_1(E-E')^{m/2} ,\\
& -\frac{1}{2} \int _{E'}^{E} \frac{a(\tilde{E})}{\sqrt{-Q(\tilde{E})}}d\tilde{E}
\sim \bar{C}_2 (E -E')^{1+m(E')-\tilde{m}(E')/2},
\end{align}
for non-zero constants $\bar{C}_1, \bar{C}_2$.
Hence we have $m=2+2m(E')-\tilde{m}(E')$. 
It follows from Proposition \ref{prop:AiBi} that, if $E' \in B_0$, then $\tilde{m}(E')=0$. Therefore, if $E' \in A_0$ (resp. $E' \in B_0$), then the multiplicity of zeros of the polynomial $D(E)$ at $E=E'$ is $2+2m(E')-\tilde{m}(E')$ (resp. $2+2m(E')$).
\end{proof}
As a corollary we have
\begin{cor} \label{cor:wpa}
Assume that $a(E') \neq 0$ for all $E'\in \cup _{i=0}^3 (A_i \cup B_i)$ and the equation $Q(E)=0$ does not have multiple roots.
We set
\begin{align}
& D(E)= \prod _{E' \in A_0} (E-E' ) \prod _{E' \in B_0} (E-E' )^{2} ,
\label{eq:DEc}\\
& N_k(E)= \prod _{E' \in A_{k}} (E-E' ) \prod _{E' \in B_{k}} (E-E' )^{2}, \quad (k=1,2,3). \label{eq:NkEc}
\end{align}
Then we have
\begin{equation}
\wp (\alpha )= e_k -\frac{4N_k(E)}{(\sum_{i=0}^3 l_i(l_i+1))^2D(E)}
\end{equation}
and $\deg _E N_k(E)= 1 +\deg _E D(E)$ for $k=1,2,3$.
\end{cor}

Now we consider the function $\kappa $. By combining Eq.(\ref{P1P6}) and Proposition \ref{prop:asym} we obtain the expression 
\begin{equation}
\kappa = \left(1 -\frac{2}{\sum_{i=0}^3 l_i(l_i+1)} \right) \frac{P_{\kappa }(E)}{\tilde{P}_{\kappa }(E)} \sqrt{-Q(E)}
\label{eq:kappaAB}
\end{equation}
where $P_{\kappa }(E)$ and $\tilde{P}_{\kappa }(E)$ are monic polynomials such that $\deg _E P_{\kappa } (E) =\deg _E \tilde{P}_{\kappa }(E) -g$. We investigate zeros of the polynomial $\tilde{P}_{\kappa }(E)$.
\begin{prop} \label{prop:kappa}
Assume that $a(E') \neq 0$ for all $E'\in \cup _{i=0}^3 (A_i \cup B_i)$ and the equation $Q(E)=0$ does not have multiple roots. Then we have
\begin{equation}
\kappa = \left(1 -\frac{2}{\sum_{i=0}^3 l_i(l_i+1)} \right) \frac{P_{\kappa }(E)}{\prod _{E' \in A_0} (E-E' ) \prod _{E' \in B_0} (E-E' )} \sqrt{-Q(E)}
\label{eq:kappaApBp}
\end{equation}
for some monic polynomials $P_{\kappa}(E)$ such that $\deg _E P_{\kappa } (E) =\deg _E (\prod _{E' \in A_0} (E-E' ) \prod _{E' \in B_0} (E-E' ) ) -g$.
\end{prop}

\begin{proof}
It follows from Eqs.(\ref{kap0}, \ref{kap123}), Proposition \ref{prop:AiBi} and Corollary \ref{cor:wpa} that, if $E' \not \in \cup _{i=0}^3 (A_i \cup B_i)$, then the function $\kappa $ is holomorphic in $E$ around $E=E'$. We investigate the degree of the poles at $E' \in \cup _{i=0}^3 (A_i \cup B_i)$.

We show that, if $k \in \{ 1,2,3 \}$ and $E' \in B_{k}$, then the function $\kappa $ is holomorphic in $E$ around $E=E'$. From Corollary \ref{cor:wpa} we have $\wp (\alpha ) \rightarrow e_{k}$ for some $k \in \{ 1,2,3\}$, $\alpha \rightarrow \omega _k$ and that $\zeta (\alpha )$ is holomorphic around $E=E'$. It follows from Eq.(\ref{alpEz}) that the function $\kappa $ is holomorphic in $E$ around $E=E'$.

Assume $E' \in B_{0}$. From Corollary \ref{cor:wpa} we have $\wp (\alpha ) -e_{k}=g_k(E)/(E-E')^2$ ($k=1,2,3$) for some rational function $g_k(E)$ such that $g_k(E') \neq 0 , \infty$. Since $\wp (x) \sim 1/x^2$ and $\zeta (x) \sim 1/x$ as $x \sim 0$, we have $\alpha \sim (E-E')/C'$ and $\zeta (\alpha ) \sim C' /(E-E')$ for some non-zero constants $C'$ as $E \rightarrow E'$. From Eq.(\ref{alpEz}), the function $\kappa $ has a single pole in $E$ around $E=E'$.

Next we consider the case $E' \in A_0$. From Corollary \ref{cor:wpa} we have $\wp (\alpha ) \sim (C') ^2 /(E-E')$ for some non-zero constants $C'$ as $E \rightarrow E'$. Then we have $\zeta (\alpha ) \sim C' /(E-E')^{1/2}$ as $E \rightarrow E'$. 
On the other hand we have 
\begin{equation}
\int _{E'}^{E} \frac{c(\tilde{E})}{\sqrt{-Q(\tilde{E})}}d\tilde{E} \sim \int _{E'}^{E} \bar{R} (\tilde{E}) (\tilde{E} -E' )^{-1/2} d\tilde{E} \sim \bar{C} (E -E')^{1/2}
\label{eq:asymp}
\end{equation}
for some functions $\bar{R} (E)$ and constants $\bar{C}$. 
Therefore we have $\kappa \sim -C' /(E-E')^{1/2}$ from Eq.(\ref{alpEz}).

For the case $E' \in A_{k}$ $(k \in \{1,2,3 \})$, we have $\wp (\alpha ) \rightarrow e_{k}$ and $\zeta (\alpha )$ is bounded around $E=E'$, and $\int _{E_0}^{E} \frac{c(\tilde{E})}{\sqrt{-Q(\tilde{E})}}d\tilde{E} $ is bounded around $E =E'$. Hence $\kappa $ is also bounded around $E\rightarrow E'$.

By combining with the fact that $\kappa $ is expressed as Eq.(\ref{eq:kappaAB}) and that the equation $Q(E)=0$ does not have multiple roots, the denominator of Eq.(\ref{eq:kappaApBp}) must be expressed as $\prod _{E' \in A_0} (E-E' ) \prod _{E' \in B_0} (E-E' )$. 
\end{proof}

\section{Twisted Heun and theta-twisted Heun polynomials} \label{sec:tHttH}

To investigate Eqs.(\ref{eq:wpaND}, \ref{eq:kappaApBp}) and sets $ A_i$, $B_i$ $(i=0,1,2,3)$, we introduce Heun, twisted Heun and theta-twisted Heun polynomials.

To study the set $A_i$ $(i=0,1,2,3)$, we consider Heun polynomials \cite{Ron}, which are related with quasi-exact solvability. For details see \cite[\S 5]{Tak2}.
For integers $\beta_0$, $\beta_1$, $\beta_2$, $\beta_3$ such that $-\sum_{i=0}^3 \beta_i \in 2 \mathbb{Z}_{\geq 0}$, $V_{\beta_0, \beta_1, \beta_2, \beta_3}$ is defined as the vector space spanned by functions $\big\{ \wp _1(x) ^{\beta _1} \wp_2(x)^{\beta _2} \wp_3(x)^{\beta _3}\wp(x)^n\big\} _{n=0, \dots ,-\sum_{i=0}^3 \beta_i /2 }$
Let  $\alpha _i \in \{ -l_i , l_i+1\}$
 $(i=0,1,2,3)$ and
\begin{gather}
 U_{\alpha _0, \alpha _1, \alpha _2, \alpha _3}=
\begin{cases}
V_{\alpha _0, \alpha _1, \alpha _2, \alpha _3},
& \sum_{i=0}^3 \alpha _i/2 \in \mathbb{Z}_{\leq 0} ;\\
V_{1-\alpha _0, 1-\alpha _1, 1-\alpha _2, 1-\alpha _3} ,
& \sum_{i=0}^3 \alpha _i /2\in \mathbb{Z}_{\geq 2} ;\\
\{ 0 \} ,& \mbox{otherwise}.
\end{cases} 
\end{gather}
If $l_0 +l_1 +l_2 +l_3$ is even, then the Hamiltonian $H$ (see Eq.(\ref{Ino}))
preserves the spaces
\begin{equation}
U_{-l_0,-l_1,-l_2,-l_3}, \quad U_{-l_0 ,-l_1,l_2+1 ,l_3+1}, \quad U_{-l_0,l_1+1,-l_2,l_3+1}, \quad U_{-l_0 ,l_1+1,l_2+1,-l_3},
\label{einvsp}
\end{equation}
and if $l_0 +l_1 +l_2 +l_3$ is odd, then the Hamiltonian $H$ preserves the spaces
\begin{equation}
U_{-l_0,-l_1,-l_2,l_3+1}, \quad U_{-l_0 ,-l_1,l_2+1,-l_3}, \quad U_{-l_0 ,l_1+1,-l_2,-l_3}, \quad U_{l_0+1,-l_1,-l_2,-l_3}. 
\label{oinvsp}
\end{equation}
\begin{prop} \label{prop:Ai}
Let $i \in \{0,1,2,3 \}$. Then there exists only one space $U_{\alpha _0, \alpha _1, \alpha _2, \alpha _3}$ in Eq.(\ref{einvsp}) or Eq.(\ref{oinvsp}) such that $\sum _{k=1}^3 \alpha _k \omega _k \equiv \omega _i$ (mod $2\omega_1 \Zint \oplus 2\omega_3 \Zint$).
Moveover the set $A_i$ coincides with the set of eigenvalues of the Hamiltonian on the space $U_{\alpha _0, \alpha _1, \alpha _2, \alpha _3}$. 
\end{prop}
\begin{proof}
Assume that $l_0 +l_1 +l_2 +l_3$ is even (resp. odd).
Let $U _{\alpha _0, \alpha _1, \alpha _2, \alpha _3}$ and $U_{\alpha  _0 ', \alpha _1 ', \alpha _2 ', \alpha _3 '}$ be different spaces in Eq.(\ref{einvsp}) (resp. Eq.(\ref{oinvsp})). Then it follows directly that $\sum _{k=1}^3 \alpha _k \omega _k \not \equiv \sum _{k=1}^3 \alpha _k ' \omega _k $ (mod $2\omega_1 \Zint \oplus 2\omega_3 \Zint$). Hence we obtain the first statement. It is shown in the proof of \cite[Theorem 3.2]{Tak1} that there are no common eigenvalues on the spaces $U _{\alpha _0, \alpha _1, \alpha _2, \alpha _3}$ and $U_{\alpha  _0 ', \alpha _1 ', \alpha _2 ', \alpha _3 '}$.

It is shown in \cite{Tak1} that, if $l_0 +l_1 +l_2 +l_3$ is even (resp. odd), then  the direct sum of the spaces in Eq.(\ref{einvsp}) (resp. Eq.(\ref{oinvsp})), which we denote by $\tilde{V}$, is the maximum finite-dimensional subspace of ${\mathcal F}$, and the set of eigenvalues of the Hamiltonian on the space $\tilde{V}$ coincides with the set of zeros of the polynomial $Q(E)$.
By combining with Proposition \ref{prop:AiBi}, it follows that the set $\cup _{i'=0}^3 A_{i '}$ coincides with the set of eigenvalues of the Hamiltonian on the space $\tilde{V}$.
Let $E \in A_i$. Then there exists a non-zero eigenfunction $\tilde{\Lambda }(x,E)$ written as Eq.(\ref{Lomega0}) or Eq.(\ref{Lomega}) which satisfies $\bar{\kappa }=0$ or $\kappa =0$. Then the function $\tilde{\Lambda }(x,E)$ has the same periodicity as the function $\left\{ \begin{array}{ll} 1, & (i=0), \\ \wp _i(x), & (i=1,2,3). \end{array} \right.$
Let $U_{\alpha _0, \alpha _1, \alpha _2, \alpha _3}$ be the space in Eq.(\ref{einvsp}) (resp. Eq.(\ref{oinvsp})) such that $\sum _{k=1}^3 \alpha _k \omega _k \equiv \omega _i$ (mod $2\omega_1 \Zint \oplus 2\omega_3 \Zint$). Then the function $f(x) \in U_{\alpha _0, \alpha _1, \alpha _2, \alpha _3}$ has the same periodicity as the function $\left\{ \begin{array}{ll} 1, & (i=0), \\ \wp _i(x), & (i=1,2,3). \end{array} \right.$ 
Thus an eigenfunction in $\tilde{V}$ whose eigenvalue is an element of $A_i $ belongs to the space $U_{\alpha _0, \alpha _1, \alpha _2, \alpha _3}$.
Therefore we obtain that the set $A_i$ coincides with the set of eigenvalues of the Hamiltonian on the space $U_{\alpha _0, \alpha _1, \alpha _2, \alpha _3}$.
\end{proof}

For $i \in \{0,1,2,3 \}$, we denote the monic characteristic polynomial of $H$ on the space $U_{\alpha _0, \alpha _1, \alpha _2, \alpha _3}$ in Proposition \ref{prop:Ai} by $H^{(i)}(E)$. We call them Heun polynomials.
It follows from Proposition \ref{prop:Ai} that the set $A_i$ ($i= 0,1,2,3 $) coincides with the set of zeros of the Heun polynomial $H^{(i)}(E)$.

Next we investigate the sets $B_i$ ($i=0,1,2,3$) and introduce twisted Heun polynomials.  Here we transform the Hamiltonian $H$. We set $z=\wp(x)$, $\phi (x)= \exp (\kappa x) (z-e_1)^{-l_1/2}(z-e_2)^{-l_2/2}(z-e_3)^{-l_3/2}$, and
$\widehat{H} _{\kappa }= \phi (x)^{-1} \circ H \circ\phi (x)$. Then we have
\begin{align}
 \widehat{H} _\kappa = & -4(z-e_1)(z-e_2)(z-e_3)
 \Big\{ \frac{d^2}{dz^2}+ \sum_{i=1}^3\frac{-l_i+\frac{1}{2}}{z-e_i}
 \frac{d}{dz} \Big\} \label{Hkap} \\
& +\big\{ -\left(l_0+l_1+l_2+l_3\right) \left( l_1+l_2+l_3-l_0-1\right) z\nonumber \\
& -\kappa ^2 +e_1(l_2+l_3)^2+e_2(l_1+l_3)^2+e_3(l_1+l_2)^2 \big\} \nonumber \\
& -2\kappa \sqrt{4(z-e_1)(z-e_2)(z-e_3)} \Big\{ \frac{d}{dz} - \frac{1}{2}\sum_{i=1}^3\frac{l_i}{z-e_i} \Big\}
. \nonumber
\end{align}

First we consider the set $B_p$. It follows from definition of $B_p$ and Propositions \ref{prop:AiBi}, \ref{prop:Lalp} that the condition $E \in B_p$ is equivalent to that there exists a solution to Eq.(\ref{InoEF}) written as Eq.(\ref{Lomega0}) with the condition $\bar{\kappa }\neq 0$. We consider the situation that the function written as Eq.(\ref{Lomega0}) with the condition $\bar{\kappa }\neq 0$ satisfies the differential equation Eq.(\ref{InoEF}).
We set 
\begin{equation}
\Phi (x,E) = \sum_{j=0}^{\hat{l}^{(0)}} \bar{a}_j (z-e_2)^j + \sqrt{(z-e_1)(z-e_2)(z-e_3)} \sum_{j=0}^{\check{l}^{(0)}} \bar{b}_j (z-e_2)^j,
\label{Lomega0h}
\end{equation}
$\bar{a}_j=0$ $(j>\hat{l}^{(0)}, \; j <0)$ and $\bar{b}_j$ $(j> \check{l}^{(0)}, \; j <0 )$, where $\hat{l}^{(0)}$ (resp. $\check{l}^{(0)}$) is the greatest integer that is less than or equal to $(l_0+l_1+l_2+l_3)/2$ (resp. $(l_0+l_1+l_2+l_3-3)/2$). Then Eq.(\ref{InoEF}) for $f(x)= \tilde{\Lambda} (x,E)$ in Eq.(\ref{Lomega0}) is equivalent to
\begin{equation}
(z-e_1)(z-e_2)(z-e_3) \left( \widehat{H} _{\bar{\kappa }}-E \right) \Phi (x,E) 
=0.
\label{eq:Hkappa}
\end{equation}
By substituting Eq.(\ref{Lomega0h}) into Eq.(\ref{eq:Hkappa}) and comparing coefficients of $(z-e_2) ^{j+4}$ and $\sqrt{(z-e_1)(z-e_2)(z-e_3)} (z-e_2)^{j+4}$, we obtain relations
\begin{align}
& (l_1+l_2+l_3-l_0-1-2j)(l_1+l_2+l_3+l_0 -2j) \bar{a}_{j} \label{eq:aj4} \\
& +\sum_{j'=j+1}^{j+4} \bar{a}_{j'} c_{\bar{a}\bar{a}}(j,j',{\bar{\kappa }}) +{\bar{\kappa }}\sum_{j'=j-1}^{j+3} \bar{b}_{j'} c_{\bar{a}\bar{b}}(j,j',{\bar{\kappa }})=0 \nonumber \\
& (l_1+l_2+l_3-l_0-4-2j)(l_1+l_2+l_3+l_0-3-2j) \bar{b}_{j} \label{eq:bj4} \\
& +\sum_{j'=j+1}^{j+4} \bar{b}_{j'} c_{\bar{b}\bar{b}}(j,j',{\bar{\kappa }}) +{\bar{\kappa }}\sum_{j'=j+2}^{j+4} \bar{a}_{j'} c_{\bar{b}\bar{a}}(j,j',{\bar{\kappa }})=0 \nonumber, 
\end{align}
where $ c_{\bar{a}\bar{a}}(j,j',{\bar{\kappa }})$, $ c_{\bar{a}\bar{b}}(j,j',{\bar{\kappa }})$, $ c_{\bar{b}\bar{a}}(j,j',{\bar{\kappa }})$ and $ c_{\bar{b}\bar{b}}(j,j',{\bar{\kappa }})$ are polynomials in $e_1$, $e_2$, $e_3$, ${\bar{\kappa }}$ and $E$.
It follows from Proposition \ref{prop:AiBi} that, if $\Phi (x,E)$ satisfies Eq.(\ref{eq:Hkappa}), ${\bar{\kappa }}\neq 0$ and $l_0+l_1+l_2+l_3$ is even (resp. odd), then the value $\bar{a}_{\hat{l}^{(0)}}$ (resp. $\bar{b}_{\check{l}^{(0)}}$) is non-zero. We set $\bar{a}_{\hat{l}^{(0)}} =1$ (resp. $\bar{b}_{\check{l}^{(0)}}=1$) for normalization.

If $l_0+l_1+l_2 +l_3$ is even (resp. odd) and $l_0> l_1+l_2+l_3-4$ (resp. $l_0 >l_1+l_2+l_3-1$), then values $\bar{b}_{\hat{l}^{(0)}-2}, \bar{a}_{\hat{l}^{(0)}-1}, \bar{b}_{\hat{l}^{(0)}-3}, \dots ,\bar{a}_{0}$ (resp. $\bar{a}_{\check{l}^{(0)}+1}, \bar{b}_{\check{l}^{(0)}-1}, \bar{a}_{\check{l}^{(0)}}, \dots ,\bar{a}_{0}$) are determined recursively by Eqs.(\ref{eq:aj4}, \ref{eq:bj4}). To satisfy Eq.(\ref{eq:Hkappa}) the coefficients of $(z-e_2)^j$ and  $\sqrt{(z-e_1)(z-e_2)(z-e_3)}  (z-e_2)^j$ $(j=0,1,2,3)$ on the l.h.s. should be zero.
Hence we have simultaneous equations for ${\bar{\kappa }}, E$, which are Eqs.(\ref{eq:aj4}, \ref{eq:bj4}) for the cases $j=-1,-2,-3,-4$ with $\bar{a}_{-4}=\bar{b}_{-4}= \dots = \bar{a}_{-1}=\bar{b}_{-1}=0$. We solve them by throwing away the solutions corresponding to the case ${\bar{\kappa }}=0$ and removing the ${\bar{\kappa }}$ factors by the elimination method in the theory of Gr\"obner bases. Then we obtain an equation $Ht^{(p)}(E)=0$, where $Ht^{(p)}(E)$ is a polynomial which is normalized to be monic.  We call $Ht^{(p)}(E)$ the twisted Heun polynomial. (See \cite{Mai} for twisted Lam\'e polynomial.) Note that we can remove completely the term including ${\bar{\kappa }}$, because it follows from Proposition \ref{thm:paek} that the number of eigenvalues $E$ whose eigenfunction is written as Eq.(\ref{Lomega0}) is finite.

Let us consider the case $l_0+l_1+l_2 +l_3$ is even (resp. odd) and $l_0 \leq  l_1+l_2+l_3-4$ (resp. $l_0 \leq l_1+l_2+l_3-1$). We set $A=\bar{b}_{(l_1+l_2+l_3-l_0-4)/2}$ (resp. $A=\bar{a}_{(l_1+l_2+l_3-l_0-1)/2}$). Then values $\bar{a}_{j}$ and $\bar{b}_{j}$ except for $\bar{b}_{(l_1+l_2+l_3-l_0-4)/2}$ (resp. $\bar{a}_{(l_1+l_2+l_3-l_0-1)/2}$) are determined recursively by Eqs.(\ref{eq:aj4}, \ref{eq:bj4}). To satisfy Eq.(\ref{eq:Hkappa}) the coefficients of $(z-e_2)^j$ and $\sqrt{(z-e_1)(z-e_2)(z-e_3)} (z-e_2) ^j$ $(j=0,1,2,3)$ on the l.h.s. should be zero.
Hence we have simultaneous equations for ${\bar{\kappa }}, E,A $. We solve them by throwing away the solutions corresponding to the case ${\bar{\kappa }}=0$ and removing the ${\bar{\kappa }}$ and $A$ factors by the elimination method in the thoery of Gr\"obner bases. Then we obtain an equation $Ht^{(p)}(E)=0$, where $Ht^{(p)}(E)$ is the twisted Heun polynomial which is normalized to be monic. Note that we can also completely remove the term including ${\bar{\kappa }}$ and $A$.

Note that elements of the sets $A_p$ sometimes appear as solutions for the recursive equation for the case ${\bar{\kappa }}=0$, although some elements in $A_p$ might not appear as solutions because the condition $\bar{a}_{\hat{l}^{(0)}}\neq 0$ (or $\bar{b}_{\check{l}^{(0)}}\neq 0$) may be broken.

We consider the sets $B_{p'}$ ($p' \in \{ 0,1,2,3 \} \setminus \{ p \} $). Let $k$ be an element of $\{1,2,3\}$ such that $\omega _{p'} \equiv \omega _{p}+ \omega _{k}$ (mod $2\omega_1 \Zint \oplus 2\omega_3 \Zint$). It follows from the definition of $B_{p'}$ and Propositions \ref{prop:AiBi}, \ref{prop:Lalp} that the condition $E \in B_{p'}$ is equivalent to that there exists a solution to Eq.(\ref{InoEF}) written as Eq.(\ref{Lomega}) with the condition $\kappa \neq 0$. We consider the situation that the function written as Eq.(\ref{Lomega}) with the condition $\kappa \neq 0$ satisfies the differential equation Eq.(\ref{InoEF}).
We set 
\begin{equation}
\Phi (x,E) = \sqrt{z-e_k} \left( \sum_{j=0}^{\hat{l}} a_j (z-e_k)^j \right) + \sqrt{(z-e_{k'})(z-e_{k''})}\left( \sum_{j=0}^{\check{l}} b_j (z-e_k)^j \right) ,
\label{Lomegah123}
\end{equation}
$a_j=0$ $(j>\hat{l}, \; j<0)$ and $b_j$ $(j> \check{l}, \; j<0 )$ where $k', k''$ are determined as $\{k,k',k''\}=\{1,2,3\}$, $z=\wp (x)$ and  $\hat{l}$ (resp. $\check{l}$) is the greatest integer that is less than or equal to $(l_0+l_1+l_2+l_3-1)/2$ (resp. $(l_0+l_1+l_2+l_3-2)/2$).
Then Eq.(\ref{InoEF}) for $f(x)=\tilde{\Lambda} (x,E)$ in Eq.(\ref{Lomega}) is equivalent to
\begin{equation}
(z-e_1)(z-e_2)(z-e_3) \left( \widehat{H} _\kappa -E \right) \Phi (x,E) 
=0,
\label{eq:Hkappa123}
\end{equation}
(see Eq.(\ref{Hkap})). By substituting Eq.(\ref{Lomega0h}) into Eq.(\ref{eq:Hkappa}) and comparing coefficients of $ \sqrt{z-e_k} (z-e_k) ^{j+4}$ and $ \sqrt{(z-e_{k'})(z-e_{k''})} (z-e_k)^{j+4}$, we obtain relations 
\begin{align}
& (l_1+l_2+l_3-l_0-2-2j)(l_1+l_2+l_3+l_0 -1-2j) a_{j} \label{eq:aj4t} \\
& +\sum_{j'=j+1}^{j+4}a_{j'} c_{aa}(j,j',\kappa ) +\kappa \sum_{j'=j}^{j+4} b_{j'} c_{ab}(j,j',\kappa )=0 \nonumber \\
& (l_1+l_2+l_3-l_0-3-2j)(l_1+l_2+l_3+l_0-2-2j) b_{j} \label{eq:bj4t} \\
& +\sum_{j'=j+1}^{j+4} b_{j'} c_{bb}(j,j',\kappa ) +\kappa \sum_{j'=j+1}^{j+3} a_{j'} c_{ba}(j,j',\kappa )=0 \nonumber, 
\end{align}
where $ c_{aa}(j,j',\kappa )$, $ c_{ab}(j,j',\kappa )$, $ c_{ba}(j,j',\kappa )$ and $ c_{bb}(j,j',\kappa )$ are polynomials in $e_1$, $e_2$, $e_3$, $\kappa$ and $E$.
It follows from Proposition \ref{prop:AiBi} that, if $\Phi (x,E)$ satisfies Eq.(\ref{eq:Hkappa123}), $\kappa \neq 0$ and $l_0+l_1+l_2+l_3$ is even (resp. odd), then the value $b_{\check{l}}$ (resp. $a_{\hat{l}}$) is non-zero. We set $b_{\check{l}} =1$ (resp. $a_{\hat{l}}=1$) for normalization.

If $l_0+l_1+l_2 +l_3$ is even (resp. odd) and $l_0> l_1+l_2+l_3-2$ (resp. $l_0 >l_1+l_2+l_3-3$), then values $a_{\check{l}-1}, b_{\check{l}-1} , a_{\check{l}-2} , \dots ,a_{0}$ (resp. $b_{\hat{l}}, a_{\hat{l}-1}, b_{\hat{l}-1}, \dots , a_{0}$) are determined recursively by Eqs.(\ref{eq:aj4t}, \ref{eq:bj4t}). To satisfy Eq.(\ref{eq:Hkappa123}) the coefficients of $ \sqrt{z-e_k} (z-e_k) ^{j}$ and $ \sqrt{(z-e_{k'})(z-e_{k''})} (z-e_k)^{j}$ $(j=0,1,2,3)$ on the l.h.s. should be zero.
Hence we have simultaneous equations for $\kappa, E$.
We solve them by throwing away the solutions corresponding to the case $\kappa=0$ and removing the $\kappa$ factors by the elimination method in the theory of Gr\"obner bases. Then we obtain an equation $Ht^{(p')}(E)=0$, where $Ht^{(p')}(E)$ is the twisted Heun polynomial which is normalized to be monic.  Note that we can also completely remove the term including $\kappa$.

If $l_0+l_1+l_2 +l_3$ is even (resp. odd) and $l_0 \leq  l_1+l_2+l_3-2$ (resp. $l_0 \leq l_1+l_2+l_3-3$), then we need to set $A=a_{(l_1+l_2+l_3-l_0-2)/2}$ (resp. $A=b_{(l_1+l_2+l_3-l_0-3)/2}$) and obtain simultaneous equation for $\kappa, E,A $. The twisted Heun polynomial $Ht^{(p')}(E)$, which is normalized to be monic, is defined similarly.

It follows from the definitions and Propositions \ref{prop:AiBi}, \ref{prop:Lalp} that the set of zeros of the Heun polynomial $H^{(i)}(E)$ coincides with the set $A_i$ and that the set of zeros of the twisted Heun polynomial $Ht^{(i)}(E)$ coincide with the set $B_i$ $(i=0,1,2,3)$.

The value $\wp (\alpha )$ is expressed by Heun and twisted Heun polynomials. 
\begin{thm} \label{thm:HHt}
Assume that the Heun polynomials $H^{(i)}(E)$ and the twisted Heun polynomials $Ht^{(i)}(E)$ do not have multiple zeros and that their zeros do not intersect with zeros of $a(E)$ for generic periods $(2\omega _1, 2\omega _3)$.
Then we have
\begin{equation}
\wp (\alpha )= e_k -\frac{4H^{(k)}(E)Ht^{(k)}(E)^2}{(\sum_{i=0}^3 l_i(l_i+1))^2H^{(0)}(E)Ht^{(0)}(E)^2}
\label{eq:thmwpa}
\end{equation}
and $\deg _E H^{(k)}(E)Ht^{(k)}(E)^2= 1 +\deg _E H^{(0)}(E)Ht^{(0)}(E)^2$ for $k=1,2,3$.
\end{thm}
\begin{proof}
Assume that $(2\omega _1, 2\omega _3)$ are the basic periods that satisfy the assumption of the theorem. Since the set $A_i$ (resp. $B_i$) ($i=0,1,2,3$) coincides with the set of zeros of the polynomial $H ^{(i)} (E)$ (resp. $Ht ^{(i)} (E)$), we obtain
\begin{equation}
 H ^{(i)} (E) =\prod _{E' \in A_i } (E-E'), \quad  Ht ^{(i)} (E) =\prod _{E' \in B_i } (E-E').
\label{eq:ABHHt}
\end{equation}
From Corollary \ref{cor:wpa}, we obtain Eq.(\ref{eq:thmwpa}) for the case that satisfies the assumption of the theorem for basic periods. On the other hand $\wp (\alpha )$ is a rational function in $E$. Hence Eq.(\ref{eq:thmwpa}) is true without generic conditions.
\end{proof}
In section \ref{sec:exa} it is shown that, if $g\leq 3$, then the assumption of Theorem \ref{thm:HHt} is true.

Now we introduce the theta-twisted Heun polynomial $H\theta (E)$ to calculate the polynomial $P_{\kappa }(E)$ appeared in Eq.(\ref{eq:kappaApBp}).

We consider the case $\sum _{i=0}^3 l_i \omega _i \equiv 0$ (mod $2\omega_1 \Zint \oplus 2\omega_3 \Zint$).
We set
\begin{align}
& u_j =\frac{\Phi _0 (x,\alpha ) (\wp (x) -e_2) ^j }{(\wp (x) -e_1)^{l_1/2}(\wp (x) -e_2)^{l_2/2}(\wp (x) -e_3)^{l_3/2}}, \\ 
& v_j = \frac{\left( \frac{\partial }{\partial x } \Phi _0 (x,\alpha ) \right) (\wp (x) -e_2 )^j}{(\wp (x) -e_1)^{l_1/2}(\wp (x) -e_2)^{l_2/2}(\wp (x) -e_3)^{l_3/2}},
\end{align}
(see Eq.(\ref{Phii})). Let $\tilde{\Lambda }(x,E)$ be the function written as Eq.(\ref{Lalpha}) and assume that $\kappa =0$ and $\alpha \not \equiv 0$ (mod $\omega_1 \Zint \oplus \omega_3 \Zint$).
It follows from periodicities and position of poles that, if $\sum _{i=0}^3 l_i \omega _i \equiv 0$ (mod $2\omega_1 \Zint \oplus 2\omega_3 \Zint$), $\kappa =0$ and $\alpha \not \equiv 0$ (mod $\omega_1 \Zint \oplus \omega_3 \Zint$), then the function $\tilde{\Lambda }(x,E)$ is written as
\begin{equation}
\tilde{\Lambda }(x,E) = \sum _{j=0}^{\hat{l}} c_j u_j + \sum _{j=0}^{\check{l}} d_j v_j
\label{eq:thtL0}
\end{equation}
where $\hat{l}$ (resp. $\check{l}$) is the greatest integer that is less than or equal to $(l_0+l_1+l_2+l_3-1)/2$ (resp. $(l_0+l_1+l_2+l_3-2)/2$).
Conversely it is shown similarly to Proposition \ref{prop:Lalp} that, if $\alpha \not \equiv 0$  (mod $\omega_1 \Zint \oplus \omega_3 \Zint$) and the r.h.s. of Eq.(\ref{eq:thtL0}) satisfies the differential equation (\ref{InoEF}), then we have $\kappa =0$.

We substitute Eq.(\ref{eq:thtL0}) into Eq.(\ref{InoEF}). By using a relation $\left( \frac{\partial }{\partial x}\right) ^2 \Phi _0 (x,\alpha ) = (2\wp (x) +\wp (\alpha )) \Phi _0 (x,\alpha )$ and comparing coefficients of $u_{j+4}$ and $v_{j+4}$, we obtain relations
\begin{align}
& (l_1+l_2+l_3-l_0-2-2j)(l_1+l_2+l_3+l_0-1-2j) c_{j} \label{eq:aj4tt} \\
& +\sum_{j'=j+1}^{j+4} c_{j'} c_{cc}(j,j',\alpha ) + \wp '(\alpha ) \sum_{j'=j+2}^{j+4} d_{j'} c_{cd}(j,j',\alpha )=0 \nonumber \\
& (l_1+l_2+l_3-l_0-3-2j)(l_1+l_2+l_3+l_0-2-2j) d_{j} \label{eq:bj4tt} \\
& +\sum_{j'=j+1}^{j+4} d_{j'} c_{dd}(j,j',\alpha ) + \wp '(\alpha ) \sum_{j'=j+1}^{j+4} c_{j'} c_{dc}(j,j',\alpha )=0 \nonumber, 
\end{align}
where $ c_{cc}(j,j',\alpha )$, $ c_{cd}(j,j',\alpha )$, $ c_{dc}(j,j',\alpha )$ and $ c_{dd}(j,j',\alpha )$ are polynomials in $e_1$, $e_2$, $e_3$, $\wp (\alpha ), \wp '(\alpha )$ and $E$.
If $\tilde{\Lambda }(x,E)$ satisfies Eq.(\ref{InoEF}), $\alpha \not \equiv 0$ (mod $\omega_1 \Zint \oplus \omega_3 \Zint$) and $l_0+l_1+l_2+l_3$ is even (resp. odd), then it follows that $Q(E)\neq 0$ and the value $d_{\check{l}}$ (resp. $c_{\hat{l}}$) is non-zero, which is obtained similarly to Proposition \ref{prop:AiBi}. We set $d_{\check{l}} =1$ (resp. $c_{\hat{l}}=1$) for normalization.

If $l_0+l_1+l_2 +l_3$ is even (resp. odd) and $l_0> l_1+l_2+l_3-2$ (resp. $l_0 >l_1+l_2+l_3-3$), then values $c_{\check{l}}, d_{\check{l}-1}, c_{\check{l}-1} , \dots , c_{0}$ (resp. $d_{\hat{l}-1}, c_{\hat{l}-1}, d_{\hat{l}-2} , \dots ,c_{0}$) are determined recursively by Eqs.(\ref{eq:aj4tt}, \ref{eq:bj4tt}). To satisfy Eq.(\ref{InoEF}) the coefficients of $u_j$ and  $v_j$ $(j=0,1,2,3)$ on the l.h.s. should be zero.
Hence we have simultaneous equations for $\wp(\alpha ), \wp '(\alpha) , E$. 
Note that the relation $\wp '(\alpha) ^2 =4(\wp(\alpha ) -e_1)(\wp(\alpha ) -e_2)(\wp(\alpha ) -e_3)$ is also satisfied.
We solve them by throwing away the solutions corresponding to the case $\alpha \equiv 0$ (mod $\omega_1 \Zint \oplus \omega_3 \Zint$) and removing the $\wp(\alpha ), \wp '(\alpha)$ factors by the elimination method in the thoery of Gr\"obner bases. Then  we obtain an equation $H\theta (E)=0$, where $H\theta (E)$ is a polynomial normalized to be monic.  We call $H \theta (E)$ the theta-twisted Heun polynomial. (see \cite{Mai} for theta-twisted Lam\'e polynomial.) Note that we can remove completely the term including $\wp (\alpha ), \wp '(\alpha )$, because it follows from Eq.(\ref{P1P6}) that the number of eigenvalues $E$ whose eigenfunction is written as a linear combination of $u_j$ and $v_j$ ($j=0,1,\dots $) is finite.

Let us consider the case $l_0+l_1+l_2 +l_3$ is even (resp. odd) and $l_0 \leq  l_1+l_2+l_3-2$ (resp. $l_0 \leq l_1+l_2+l_3-3$). We set $A=c_{(l_1+l_2+l_3-l_0-2)/2}$ (resp. $A=d_{(l_1+l_2+l_3-l_0-3)/2}$). Then values $c_{j}$ and $d_{j}$ except for $c_{(l_1+l_2+l_3-l_0-2)/2}$ (resp. $d_{(l_1+l_2+l_3-l_0-3)/2}$) are determined recursively by Eqs.(\ref{eq:aj4tt}, \ref{eq:bj4tt}). To satisfy Eq.(\ref{InoEF}) the coefficients of $u_j$ and $v_j$ $(j=0,1,2,3)$ on the l.h.s. should be zero.
Hence we have simultaneous equations for $\wp (\alpha ), \wp '(\alpha ), E,A $. We solve them by throwing away the solutions corresponding to the case $\alpha \equiv 0$ (mod $\omega_1 \Zint \oplus \omega_3 \Zint$) and removing the $\wp(\alpha ), \wp '(\alpha), A$ factors by the elimination method in the theory of Gr\"obner bases. Then we obtain an equation $H\theta (E)=0$, where $H\theta (E)$ is the theta-twisted Heun polynomial which is normalized to be monic.
Note that we can also remove completely the term including $\wp (\alpha ), \wp '(\alpha )$, $A$.

Next we consider the case $\sum _{i=0}^3 l_i \omega _i \equiv \omega _p$ (mod $2\omega_1 \Zint \oplus 2\omega_3 \Zint$) and $p\neq 0$.
We set 
\begin{equation}
\Psi _p(x,\alpha )= \frac{\sigma (x-\alpha -\omega _p)}{\sigma (x)}\exp \left( \left( \zeta (\alpha +\omega _p)  +\frac{1}{2}\frac{\wp '(\alpha +\omega _p)}{\wp (\alpha +\omega _p)-e_p} \right) x \right).
\end{equation}
Then the function $\Psi _p(x,\alpha )$ has the same periodicities as the function $\Phi _0 (x,\alpha ) (\wp (x) -e_1)^{-l_1/2}(\wp (x) -e_2)^{-l_2/2}(\wp (x) -e_3)^{-l_3/2}$, which follows from Eq.(\ref{adzeta}) for the case $z=x+\omega_p$ and $\tilde{z}=-\omega _p$. 
We set
\begin{align}
& u_j =\frac{\Psi _p(x,\alpha ) (\wp (x) -e_2) ^j }{(\wp (x) -e_1)^{l_1/2}(\wp (x) -e_2)^{l_2/2}(\wp (x) -e_3)^{l_3/2}}, \\ 
& v_j = \frac{\left( \frac{\partial }{\partial x } \Psi _p(x,\alpha ) \right) (\wp (x) -e_2 )^j}{(\wp (x) -e_1)^{l_1/2}(\wp (x) -e_2)^{l_2/2}(\wp (x) -e_3)^{l_3/2}},
\end{align}
Let $\tilde{\Lambda }(x,E)$ be the function written as Eq.(\ref{Lalpha}) and assume that $\kappa =0$ and $\alpha \not \equiv 0$ (mod $\omega_1 \Zint \oplus \omega_3 \Zint$).
It follows from periodicities and position of poles that, if $\kappa =0$ and $\alpha \not \equiv 0$ (mod $\omega_1 \Zint \oplus \omega_3 \Zint$), then the function $\tilde{\Lambda }(x,E)$ is written as
\begin{equation}
\tilde{\Lambda }(x,E) = \sum _{j=0}^{\hat{l}} c_j u_j + \sum _{j=0}^{\check{l}} d_j v_j . \label{eq:thtL}
\end{equation}
Conversely it is shown similarly to Proposition \ref{prop:Lalp} that, if $\alpha \not \equiv 0$ (mod $\omega_1 \Zint \oplus \omega_3 \Zint$) and the r.h.s. of Eq.(\ref{eq:thtL}) satisfies the differential equation (\ref{InoEF}), then we have $\kappa =0$.

We substitute Eq.(\ref{eq:thtL}) into Eq.(\ref{InoEF}). By using relations of elliptic functions and comparing coefficients of $u_{j+4}$ and $v_{j+4}$, we obtain relations
\begin{align}
& ( \wp (\alpha + \omega _p) -e_p)^2 (l_1+l_2+l_3-l_0-2-2j)(l_1+l_2+l_3+l_0-1-2j) c_{j} \label{eq:aj4t123} \\
& +\sum_{j'=j+1}^{j+4} c_{j'} c_{cc}(j,j',\alpha ) + \wp '(\alpha + \omega _p) ( \wp (\alpha + \omega _p) -e_p) \sum_{j'=j+1}^{j+4} d_{j'} c_{cd}(j,j',\alpha )=0 \nonumber \\
& ( \wp (\alpha + \omega _p) -e_p)^2 (l_1+l_2+l_3-l_0-3-2j)(l_1+l_2+l_3+l_0-2-2j) d_{j-4} \label{eq:bj4t123} \\
& +\sum_{j'=j+1}^{j+4} d_{j'} c_{dd}(j,j',\alpha ) + \wp '(\alpha + \omega _p) ( \wp (\alpha + \omega _p) -e_p) \sum_{j'=j}^{j+4} c_{j'} c_{dc}(j,j',\alpha )=0 \nonumber, 
\end{align}
where $ c_{cc}(j,j',\alpha )$, $ c_{cd}(j,j',\alpha )$, $ c_{dc}(j,j',\alpha )$ and $ c_{dd}(j,j',\alpha )$ are polynomials in $e_1$, $e_2$, $e_3$, $\wp (\alpha + \omega _p ), \wp '(\alpha + \omega _p )$ and $E$.
It is shown similarly that, if $\tilde{\Lambda }(x,E)$ satisfies Eq.(\ref{InoEF}), $\alpha \not \equiv 0$ (mod $\omega_1 \Zint \oplus \omega_3 \Zint$) and $l_0+l_1+l_2+l_3$ is even (resp. odd), then the value $d_{\check{l}}$ (resp. $c_{\hat{l}}$) is non-zero. We set $d_{\check{l}} = ( \wp (\alpha + \omega _p) -e_p)^{2\check{l}} $ (resp. $c_{\hat{l}}= ( \wp (\alpha + \omega _p) -e_p)^{2\hat{l}}$) for normalization.

Similarly to the case $\sum _{i=0}^3 l_i \omega _i \equiv 0$ (mod $2\omega_1 \Zint \oplus 2\omega_3 \Zint$), the coefficients $c_j$ and $d_j$ are determined and we have simultaneous equations for $\wp (\alpha + \omega _p), \wp '(\alpha + \omega _p) ,E$ and so on to satisfy Eq.(\ref{InoEF}). Note that, if $l_0+l_1+l_2 +l_3$ is even (resp. odd) and $l_0 \leq  l_1+l_2+l_3-2$ (resp. $l_0 \leq l_1+l_2+l_3-3$), then we need to set $A=c_{(l_1+l_2+l_3-l_0-2)/2}$ (resp. $A=d_{(l_1+l_2+l_3-l_0-3)/2}$) and obtain simultaneous equations for $\wp (\alpha + \omega _p), \wp '(\alpha + \omega _p) , E,A $, else we obtain simultaneous equations for $\wp (\alpha + \omega _p), \wp '(\alpha + \omega _p) , E$. 
By throwing away the solution $\alpha \equiv 0$ (mod $\omega_1 \Zint \oplus \omega_3 \Zint$) and removing the factors except for $E$ by the elimination method in the theory of Gr\"obner bases, we obtain the theta-twisted Heun polynomial $H\theta (E)$, which is normalized to be monic and the equation $H\theta (E)=0$.

\begin{thm} \label{thm:ttH}
Assume that for generic periods $(2\omega _1, 2\omega _3)$ the polynomials $a(E)$, $Q(E)$, $H\theta (E)$, $H^{(i)} (E)$, $Ht^{(i)}(E)$, $(i=0,1,2,3)$ do not have multiple zeros, and their zeros do not intersect each other.
If  $\deg _E H\theta (E) - \deg _E H^{(0)}(E) Ht^{(0)}(E)=-g$, then we have
\begin{equation}
\kappa = \left(1 -\frac{2}{\sum_{i=0}^3 l_i(l_i+1)} \right) \frac{ H\theta (E)}{H^{(0)}(E) Ht^{(0)}(E)} \sqrt{-Q(E)}.
\label{eq:kappattH}
\end{equation}
\end{thm}
\begin{proof}
Assume that $(2\omega _1, 2\omega _3)$ are the basic periods that satisfy the assumption of the theorem. From Proposition \ref{prop:kappa} and Theorem \ref{thm:HHt}, we have
\begin{equation}
\kappa = \left(1 -\frac{2}{\sum_{i=0}^3 l_i(l_i+1)} \right) \frac{ P_{\kappa }(E)}{H^{(0)}(E) Ht^{(0)}(E)} \sqrt{-Q(E)}.
\label{eq:kappattH1}
\end{equation}
If the eigenvalue $E'$ such that $\alpha \not \equiv 0$ $($mod $\omega_1 \Zint \oplus \omega_3 \Zint)$ satisfies $H \theta (E') =0$, then we have $\kappa =0$ and $P_{\kappa }(E')=0$. By assumption, the equation $H\theta (E) =0$ does not have multiple roots. 
Hence the polynomial $P_{\kappa } (E)$ is divisible by $H\theta (E)$.
On the other hand we have $\deg _E H\theta (E) = -g+ \deg _E H^{(0)}(E) Ht^{(0)}(E) = \deg _E  P_{\kappa } (E) $ from the assumption, Proposition \ref{prop:kappa} and Eq.(\ref{eq:ABHHt}). Hence we have $P_{\kappa } (E) = H\theta (E)$ and Eq.(\ref{eq:kappattH}) for the case that satisfies the assumption for basic periods. On the other hand $\kappa /\sqrt{-Q(E)} $ is a rational function in $E$. Hence Eq.(\ref{eq:kappattH}) is true without generic conditions.
\end{proof}
In section \ref{sec:exa} we show that, if $g\leq 3$, then the assumption of Theorem \ref{thm:HHt} is true.

\section{Examples} \label{sec:exa}
In this section we explicitly calculate covering maps and several functions that have appeared in this paper. In section \ref{ssec:review} we review constructions of covering maps and related functions. In section \ref{ssec:dcc} we observe relations among different coupling constants $(l_0 , l_1 , l_2 , l_3)$. In sections \ref{ssec:g1}-\ref{ssec:g4} we express covering maps and related functions for all the cases where the genus of the hyperelliptic curve is less than or equal to three, and some cases of genus equal to four. By substituting them into the ones in section \ref{ssec:review}, we obtain several explicit formulae. Especially , for each case the monodromies are expressed as Eq.(\ref{sehypellint}) and the transformation formula between elliptic integrals of the first kind (resp. second kind) and the hyperelliptic integrals of the first kind (resp. second kind) are expressed as Eq.(\ref{alpints4}) (resp. Eq.(\ref{kap123s4}) or Eq.(\ref{kap0s4})).

\subsection{Review of covering maps and related functions} \label{ssec:review}
We review constructions of the covering maps and related functions; later we explicitly express their form for each case.

The doubly periodic function $\Xi (x,E)$ is defined in Proposition \ref{prop:prod}.
Based on the function $\Xi (x,E)$, the polynomials $Q(E)$, $c(E)$ and $a(E)$ are determined in Eqs.(\ref{const}, \ref{FFx}, \ref{polaE}).
Let $ \tilde{\Lambda }(x,E)$ be the solutions to Eq.(\ref{InoEF}) which is expressed as Eq.(\ref{eq:tilL}) and $E_0$ be an eigenvalue such that $Q(E_0)=0$.
Then there exists $q_1,q_2, q_3 \in \{0,1\}$ such that $\Lambda (x+2\omega _k,E_0)=(-1)^{q_k} \Lambda (x,E_0)$ for $k=1,2,3$.
The monodromies of solutions to Eq.(\ref{InoEF}) are expressed by hyperelliptic integrals as
\begin{equation}
\tilde{\Lambda }(x+2\omega _k,E)=(-1)^{q_k} \tilde{\Lambda }(x,E) \exp \left( -\frac{1}{2} \int_{E_0}^{E}\frac{ -2\eta _k a(\tilde{E}) +2\omega _k c(\tilde{E}) }{\sqrt{-Q(\tilde{E})}} d\tilde{E}\right) 
\label{sehypellint}
\end{equation}
for $k=1,2,3$ (see Eq.(\ref{hypellint})).
The Heun polynomials $H^{(i)}(E)$ $(i=0,1,2,3)$, the twisted Heun polynomials $Ht^{(i)}(E)$ $(i=0,1,2,3)$ and the theta-twisted Heun polynomial $H \theta (E)$ are introduced in section \ref{sec:tHttH}. 
Let $g$ be the genus of the hyperelliptic curve $\nu ^2 = -Q(E)$.
From explicit expressions of functions $a(E)$, $c(E)$, $H\theta (E)$, $H^{(i)}(E)$ and $Ht^{(i)} (E)$ $(i=0,1,2,3)$ given in a list in sections \ref{ssec:g1}-\ref{ssec:g3}, assumptions of Theorem \ref{thm:HHt} and Theorem \ref{thm:ttH} are shown to be correct for the case $g \leq 3$. 
It can be shown similarly that they are also correct for the cases $l_0 \leq 10$ and $l_1=l_2=l_3=0$.
In the case where assumptions of Theorem \ref{thm:HHt} and Theorem \ref{thm:ttH} are correct, the image of the covering map $\xi (=\wp (\alpha ))$ and $\kappa$ are calculated as in Eqs.(\ref{eq:thmwpa}, \ref{eq:kappattH}), i.e.,
\begin{align}
& \xi  = e_k -\frac{4H^{(k)}(E)Ht^{(k)}(E)^2}{(\sum_{i=0}^3 l_i(l_i+1))^2H^{(0)}(E)Ht^{(0)}(E)^2} \quad (k=1,2,3),\label{wpas4} \\
& \kappa = \left(1 -\frac{2}{\sum_{i=0}^3 l_i(l_i+1)} \right) \frac{ H\theta (E)}{H^{(0)}(E) Ht^{(0)}(E)} \sqrt{-Q(E)}. \label{kappas4}
\end{align}
Then a formula which reduces a hyperelliptic integral of the first kind to an elliptic integral of the first kind, i.e.,
\begin{equation}
\int _{\infty} ^{\xi } \frac{d\tilde{\xi } }{\sqrt{4\tilde{\xi } ^3-g_2\tilde{\xi } -g_3}} = -\frac{1}{2} \int _{\infty}^{E} \frac{a(\tilde{E})}{\sqrt{-Q(\tilde{E})}}d\tilde{E} ( = \alpha )
\label{alpints4}
\end{equation}
is obtained (see Eq.(\ref{alpint})). Next we express a formula which reduces a hyperelliptic integral of the second kind to an elliptic integral of the second kind. Since the value  $E_0$ satisfies $Q(E_0)=0$ we have $E_0 \in A_i$ for some $i\in \{ 0,1,2,3 \}$ (see Proposition \ref{prop:AiBi}).
Let $\xi _0 $ be the value $\xi $ evaluated at $E=E_0$ determined by Eq.(\ref{wpas4}). 
If $E_0 \in A_1 \cup A_2 \cup A_3$, then we also obtain a formula which reduces a hyperelliptic integral of the second kind to an elliptic integral of the second kind, i.e.,
\begin{equation}
\frac{1}{2} \int  _{E_0}^{E} \frac{c(\tilde{E})}{\sqrt{-Q(\tilde{E})}}d\tilde{E} = -\kappa + \int _{\xi _0 } ^{\xi } \frac{\tilde{\xi } d\tilde{\xi } }{\sqrt{4\tilde{\xi } ^3-g_2\tilde{\xi } -g_3}} . \label{kap123s4}
\end{equation}
(see Eq.(\ref{kap123})). If $E_0 \in A_0$, then we have 
\begin{equation}
\kappa = -\frac{1}{2} \int  _{E_0}^{E} \frac{c(\tilde{E})}{\sqrt{-Q(\tilde{E})}}d\tilde{E} + \zeta \left(  \frac{1}{2}\int_{E_0}^{E}\frac{a(\tilde{E})}{\sqrt{-Q(\tilde{E})}} d\tilde{E} \right) . \label{kap0s4}
\end{equation}
(see Eq.(\ref{kap0})). Thus the solution to Eq.(\ref{InoEF}) is expressed in the form of the Hermite-Krichever Ansatz (see Eq.(\ref{Lalpha}))
and the monodromies are expressed as
\begin{align} 
& \tilde{\Lambda }(x+2\omega _k,E)  = \exp (-2\eta _k \alpha +2\omega _k \zeta (\alpha ) +2 \kappa \omega _k ) \tilde{\Lambda }(x,E) \label{ellints4}
\end{align}
for $k=1,2,3$ (see Eq.(\ref{ellint})), where the value $\alpha $ is determined from Eq.(\ref{alpints4}). The value $\kappa $ is written as Eq.(\ref{kappas4}).
By expressing $\alpha $ and $\zeta (\alpha )$ as elliptic integrals in variable $\xi (=\wp (\alpha ))$, the monodromies are represented in terms of elliptic integrals through a transformation by the covering map $\xi$ (see Eq.(\ref{wpas4})).

\subsection{Relations among the cases of different coupling constants} \label{ssec:dcc}

We comment on relationship among several coupling constants $(l_0, l_1, l_2, l_3)$.
Let $\Xi _{l_0, l_1, l_2, l_3} ^{( 2\omega _1, 2\omega _3)} (x,E)$ (resp. $Q _{l_0, l_1, l_2, l_3} ^{( 2\omega _1, 2\omega _3)}(x)$, $a _{l_0, l_1, l_2, l_3} ^{( 2\omega _1, 2\omega _3)}(x)$, $c _{l_0, l_1, l_2, l_3} ^{( 2\omega _1, 2\omega _3)}(x)$, $\xi _{l_0, l_1, l_2, l_3} ^{( 2\omega _1, 2\omega _3)}$, $\kappa _{l_0, l_1, l_2, l_3} ^{( 2\omega _1, 2\omega _3)} $) be the functions $\Xi (x,E)$ (resp. $Q(x)$, $a(x)$, $c(x)$, $\xi$, $\kappa $) for the periods $(2\omega_1, 2\omega _3)$ and the coupling constants $(l_0, l_1, l_2, l_3)$.

The following proposition is obtained by parallel transformation $x \rightarrow x+\omega _k$ $(k=1,2,3)$ and definitions of functions (see Eqs.(\ref{const}, \ref{FFx}, \ref{polaE}, \ref{alpints4}, \ref{kap123s4}, \ref{kap0s4})):
\begin{prop} \label{prop:ptrans}
We have
\begin{align}
& \Xi _{l_1, l_0, l_3, l_2} ^{( 2\omega _1, 2\omega _3)} (x ,E) = \Xi _{l_0, l_1, l_2, l_3} ^{( 2\omega _1, 2\omega _3)} (x +\omega _1,E), \quad \Xi _{l_2, l_3, l_0, l_1} ^{( 2\omega _1, 2\omega _3)} (x ,E) = \Xi _{l_0, l_1, l_2, l_3} ^{( 2\omega _1, 2\omega _3)} (x +\omega _2,E) ,\\
& \Xi _{l_3, l_2, l_1, l_0} ^{( 2\omega _1, 2\omega _3)} (x ,E) = \Xi _{l_0, l_1, l_2, l_3} ^{( 2\omega _1, 2\omega _3)} (x +\omega _3,E) ,\nonumber \\
& X _{l_0, l_1, l_2, l_3} ^{( 2\omega _1, 2\omega _3)}(x) = X _{l_1, l_0, l_3, l_2} ^{( 2\omega _1, 2\omega _3)}(x) = X _{l_2, l_3, l_0, l_1} ^{( 2\omega _1, 2\omega _3)}(x) = X _{l_3, l_2, l_1, l_0} ^{( 2\omega _1, 2\omega _3)}(x), \nonumber \\
& \chi _{l_0, l_1, l_2, l_3} ^{( 2\omega _1, 2\omega _3)} = \chi  _{l_1, l_0, l_3, l_2} ^{( 2\omega _1, 2\omega _3)} = \chi _{l_2, l_3, l_0, l_1} ^{( 2\omega _1, 2\omega _3)} = \chi _{l_3, l_2, l_1, l_0} ^{( 2\omega _1, 2\omega _3)}, \nonumber 
\end{align}
for $X= Q$, $a$, $c$ and $\chi = \xi $, $\kappa $.
\end{prop}
By Propotision \ref{prop:ptrans} and changing periods, it follows that
\begin{prop} \label{prop:pchange}
We have
\begin{align}
& \Xi _{l_0, l_1, l_3, l_2} ^{( 2\omega _1, -2\omega _1-2\omega _3)} (x ,E) = \Xi _{l_1, l_0, l_2, l_3} ^{( 2\omega _1, -2\omega _1-2\omega _3)} (x +\omega _1 ,E) = \Xi _{l_0, l_1, l_2, l_3} ^{( 2\omega _1, 2\omega _3)} (x,E), \\
& \Xi _{l_0, l_3, l_2, l_1} ^{( 2\omega _3, 2\omega _1)} (x ,E) = \Xi _{l_2, l_1, l_0, l_3} ^{( 2\omega _3, 2\omega _1)} (x +\omega _2 ,E) = \Xi _{l_0, l_1, l_2, l_3} ^{( 2\omega _1, 2\omega _3)} (x,E), \nonumber \\
& \Xi _{l_0, l_2, l_1, l_3} ^{( -2\omega _1-2\omega _3 , 2\omega _3)} (x ,E) = \Xi _{l_3, l_1, l_2, l_0} ^{( -2\omega _1-2\omega _3, 2\omega _3)} (x +\omega _3 ,E) = \Xi _{l_0, l_1, l_2, l_3} ^{( 2\omega _1, 2\omega _3)} (x,E), \nonumber \\
& X _{l_0, l_1, l_3, l_2} ^{( 2\omega _1, -2\omega _1-2\omega _3)} (x) = X _{l_1, l_0, l_2, l_3} ^{( 2\omega _1, -2\omega _1-2\omega _3)} (x) = X _{l_0, l_1, l_2, l_3} ^{( 2\omega _1, 2\omega _3)} (x),  \nonumber \\
& X _{l_0, l_3, l_2, l_1} ^{( 2\omega _3, 2\omega _1)} (x) = X _{l_2, l_1, l_0, l_3} ^{( 2\omega _3, 2\omega _1)} (x) = X _{l_0, l_1, l_2, l_3} ^{( 2\omega _1, 2\omega _3)} (x) ,  \nonumber \\
& X _{l_0, l_2, l_1, l_3} ^{( -2\omega _1-2\omega _3 , 2\omega _3)} (x) = X _{l_3, l_1, l_2, l_0} ^{( -2\omega _1-2\omega _3, 2\omega _3)} (x ) = X _{l_0, l_1, l_2, l_3} ^{( 2\omega _1, 2\omega _3)} (x),  \nonumber \\
& \chi _{l_0, l_1, l_3, l_2} ^{( 2\omega _1, -2\omega _1-2\omega _3)}  = \chi _{l_1, l_0, l_2, l_3} ^{( 2\omega _1, -2\omega _1-2\omega _3)} = \chi _{l_0, l_1, l_2, l_3} ^{( 2\omega _1, 2\omega _3)} ,  \nonumber \\
&  \chi  _{l_0, l_3, l_2, l_1} ^{( 2\omega _3, 2\omega _1)} = \chi  _{l_2, l_1, l_0, l_3} ^{( 2\omega _3, 2\omega _1)} = \chi  _{l_0, l_1, l_2, l_3} ^{( 2\omega _1, 2\omega _3)} , \nonumber \\
&  \chi  _{l_0, l_2, l_1, l_3} ^{( -2\omega _1-2\omega _3 , 2\omega _3)}  = \chi  _{l_3, l_1, l_2, l_0} ^{( -2\omega _1-2\omega _3, 2\omega _3)} = \chi  _{l_0, l_1, l_2, l_3} ^{( 2\omega _1, 2\omega _3)} , \nonumber 
\end{align}
for $X= Q$, $a$, $c$ and $\chi = \xi $, $\kappa $.
\end{prop}
We can replace the coupling constants $(l_0, l_1, l_2, l_3)$ to the case $l_0 \geq l_1 \geq l_2 \geq l_3$ by repeatedly applying Proposition \ref{prop:pchange}.
We consider examples mainly for the case $l_0 \geq l_1 \geq l_2 \geq l_3$ in sections \ref{ssec:g1}-\ref{ssec:g4}.

The relationship between the case $(l_0 ,0,0,0)$ and the case $( l_0 ,l_0 , l_0 ,l_0)$ is connected by changing the scale of the periods $(2\omega _1 , 2\omega _3) \leftrightarrow (\omega _1 , \omega _3)$, and the relationship between the cases $(l_0 ,l_1,0,0)$ and $( l_0 ,l_0 ,l_1 ,l_1)$ is connected by changing periods $(2\omega _1 , 2\omega _3) \leftrightarrow (\omega _3 , 2\omega _1)$.

\subsection{The case $g=1$} \label{ssec:g1}
 
If the genus of the related hyperelliptic curve $\nu ^2=-Q(E)$ is one and the condition $l_0 \geq l_1 \geq l_2 \geq l_3$ is satisfied, then we have three cases 
$$
(l_0, l_1 , l_2, l_3) =(1,0,0,0), \; (1,1,0,0), \; (1,1,1,1).
$$
\subsubsection{The case $(l_0, l_1 , l_2, l_3) =(1,0,0,0)$}
We have \\
$ \; \; \Xi (x,E)=E+\wp(x), \quad a(E)=1, \quad c(E)=E , \quad  Q(E)=(E+e_1)(E+e_2)(E+e_3), $\\
$ \; \; H^{(0)}(E)=1, \quad H^{(k)} (E)=E+e_k, \; (k=1,2,3), $\\
$ \; \; Ht ^{(i)} (E)=1 , \; (i=0,1,2,3), \quad H\theta (E)=1.$\\
For this case, an eigenfunction of the Hamiltonian $H$ with eigenvalue $E$ is expressed as
\begin{equation*}
\exp(\zeta (\alpha )x) \sigma (x-\alpha )/\sigma (x ), \; \; E= -\wp(\alpha ),
\end{equation*}
and we have $\kappa =1$ and $\xi=-E$. The formulae (\ref{alpints4}, \ref{kap123s4}) are trivial, because they only change signs by $\xi =-E$.

\subsubsection{The case $(l_0, l_1 , l_2, l_3) =(1,1,0,0)$}
We have \\
$ \; \; \Xi (x,E)=E+\wp(x) +\wp( x +\omega_1) -3e_1 , \quad a(E)=2, \quad c(E)=E-3e_1, $\\
$ \; \; H^{(0)} (E)=E-4e_1 , \quad H^{(1)}(E)=E^2-2e_1E+g_2-11e_1^2 , \quad  H^{(2)} (E)= H^{(3)} (E)=1, $\\
$ \; \; Q(E)= H^{(0)}(E) H^{(1)}(E), \quad   Ht ^{(0)} (E)= Ht ^{(1)} (E)= 1 ,$\\
$ \; \; Ht ^{(2)} (E)=E+5e_2+3e_3, \quad  Ht ^{(3)} (E)=E+3e_2+5e_3, \quad H\theta (E)= 1.$\\
The formulae (\ref{alpints4}, \ref{kap0s4}) are related with Landen transformation, i.e. changing periods $(2\omega _1, 2\omega _3 ) \leftrightarrow (\omega _1, 2\omega _3 )$.

\subsubsection{The case $(l_0, l_1 , l_2, l_3) =(1,1,1,1)$}
We have \\
$ \; \; \Xi (x,E)=E+\wp(x) +\wp (x+\omega _1)+\wp (x+\omega _2) +\wp (x+\omega _3), $\\
$ \; \; a(E)=4, \quad c(E)=E, \quad  Q(E)=(E+4e_1)(E+4e_2)(E+4e_3), $\\
$ \; \; H^{(0)}(E)=Q(E), \quad H^{(k)} (E)=1, \; (k=1,2,3), \quad Ht ^{(0)} (E)=1 ,$\\
$ \; \; Ht ^{(k)} (E)= E^2+8e_kE+4g_2-32e_k^2, \; (k=1,2,3), \quad H\theta (E)= E^2 -\frac{4}{3}g_2.$\\
The formulae (\ref{alpints4}, \ref{kap0s4}) are related with the changing scales of the periods $(2\omega _1, 2\omega _3 ) \leftrightarrow (\omega _1, \omega _3 )$.

\subsection{The case $g=2$} \label{ssec:g2}
If the genus of the related hyperelliptic curve $\nu ^2=-Q(E)$ is two and the condition $l_0 \geq l_1 \geq l_2 \geq l_3$ is satisfied, then we have seven cases 
$$
(l_0, l_1 , l_2, l_3) = \left\{ 
\begin{array}{l}
(2,0,0,0), \; (2,1,0,0), \; (2,1,1,0), \; (2,2,0,0), \\
(2,2,1,1), \; (2,2,2,2), \; (1,1,1,0).
\end{array} \right.
$$

\subsubsection{The case $(l_0, l_1 , l_2, l_3) =(2,0,0,0)$} $ $\\
$ \; \; \Xi (x,E)=9\wp(x)^2+3E\wp(x)+E^2-\frac{9}{4}g_2, \quad a(E)=3E, \quad c(E)=E^2-\frac{3}{2}g_2, $\\
$ \; \; \textstyle H^{(0)}(E)=E^2-3g_2, \quad H^{(k)} (E)=E-3e_k, \; (k=1,2,3), \quad  Q(E)= \prod_{i=0}^3 H^{(i)}(E)$\\
$ \; \; Ht ^{(0)} (E)=1 , \quad Ht ^{(k)} (E)=E+6e_k , \; (k=1,2,3), \quad H\theta (E)=1$.

By setting $\xi =-y/6$, $E=z$, $g_2=a/3$, $g_3=b/54$, Eq.(\ref{Hermitef}) is obtained from Eq.(\ref{alpints4}), and Eq.(\ref{Hermites}) is obtained from Eq.(\ref{kap123s4}).

\subsubsection{The case $(l_0, l_1 , l_2, l_3) =(2,1,0,0)$} $ $\\
$ \; \;  \Xi (x,E)=9\wp(x)^2+3(E-2e_1)\wp(x)+(E+4e_1)\wp (x+\omega _1)+E^2-5e_1E-27e_1^2 , $\\
$ \; \;  a(E)=4E-2e_1, \quad c(E)=E^2-5e_1E-27e_1^2+\frac{3}{4}g_2 , \quad  Q(E)= \prod_{i=0}^3 H^{(i)}(E),$\\
$ \; \; H^{(0)} (E)=E+4e_1, \quad  H^{(2)} (E)=E^2+2(4e_3+3e_2)E-52e_2e_3-12e_3^2-31e_2^2,$\\
$ \; \;  H^{(1)}(E)=1, \quad H^{(3)} (E)=E^2+2(4e_2+3e_3)E-52e_2e_3-12e_2^2-31e_3^2,$ \\
$ \; \; Ht ^{(0)} (E)=E-\frac{19}{2}e_1  , \quad Ht ^{(1)} (E)=E^2-e_1E-101e_1^2+\frac{27}{4}g_2,$\\
$ \; \;  Ht ^{(2)} (E)=E+14e_2+5e_3, \quad Ht ^{(3)} (E)=E+5e_2+14e_3 , \quad H\theta (E)=1.$

\subsubsection{The case $(l_0, l_1 , l_2, l_3) =(2,1,1,0)$} $ $\\
$ \begin{array}{cl} \Xi (x,E)=& 9\wp(x)^2+3(E+2e_3)\wp(x)+ (E+e_1-8e_2) \wp (x+\omega_1) \\
& + (E+e_2-8e_1) \wp (x+\omega_2)+ E^2+5e_3E+48e_3^2-\frac{63}{4}g_2 , 
\end{array}$\\
$ \; \;  a(E)=5E+13e_3, \quad c(E)=E^2+5e_3E+48e_3^2-15g_2,  \quad Q(E)= \prod_{i=0}^3 H^{(i)}(E), $\\
$ \; \;   H^{(0)} (E)=1  ,\quad H^{(1)} (E)= E-8e_2+e_1 , \quad  H^{(2)} (E)=E+e_2-8e_1 , $\\
$ \; \;   H^{(3)}(E)=E^3+3e_3E^2+(51e_3^2-20g_2)E+\frac{97}{4}e_3g_2-\frac{271}{4}g_3 ,$\\
$ \; \;  Ht ^{(0)} (E)=E^2+\frac{26}{5}e_3E+\frac{293}{5}e_3^2-\frac{108}{5}g_2, $\\
$ \; \;  Ht ^{(1)} (E)=E^2+2(-e_2+3e_1)E+33e_2^2+74e_1e_2-103e_1^2,$\\
$ \; \;  Ht ^{(2)} (E)=E^2+2(-e_1+3e_2)E+33e_1^2+74e_1e_2-103e_2^2 ,$\\
$ \; \;  Ht ^{(3)} (E)=E+17e_3, \quad H\theta (E)=1.$

\subsubsection{The case $(l_0, l_1 , l_2, l_3) =(2,2,0,0)$} $ $\\
$ \begin{array}{cl} \Xi (x,E)= & 9\wp(x)^2+3(E-6e_1)\wp(x)+9\wp(x+\omega_1)^2+3(E-6e_1)\wp(x+\omega _1) \\
& +E^2-15e_1E-36e_1^2+\frac{9}{2}g_2, 
\end{array}$\\
$ \; \;  a(E)=6(E-6e_1), \quad c(E)=E^2-15e_1 E-36e_1^2+6g_2, $\\
$ \; \; H^{(0)}(E)=E(E^2-12e_1E-144e_1^2+12g_2) , \quad H^{(1)}(E)=E^2-18e_1E-27e_1^2+9g_2 ,$\\
$ \; \; H^{(k)} (E)=1, \; (k=2,3), \quad  Q(E)= H^{(0)}(E)H^{(1)}(E) , \quad $\\
$ \; \; Ht ^{(0)} (E)=E-18e_1 , \quad Ht ^{(1)} (E)=E^2-432e_1^2+36g_2 , $\\
$ \; \;  Ht ^{(2)} (E) =E^3+9(5e_2+3e_3) E^2 +216e_1(e_1-e_2)E+432(e_1-e_2)(e_2 -e_3)(e_3-e_1), $\\
$ \; \;  Ht ^{(3)} (E) =E^3+9(3e_2+5e_3) E^2 +216e_1(e_1-e_3)E-432(e_1-e_2)(e_2 -e_3)(e_3-e_1), $\\
$ \; \; H\theta (E)=E^2-\frac{72}{5}e_1E-\frac{432}{5}e_1^2+\frac{36}{5}g_2. $

\subsubsection{The case $(l_0, l_1 , l_2, l_3) =(2,2,1,1)$} $ $\\
$ \begin{array}{cl} \Xi (x,E)= & 9(\wp(x)^2+\wp(x+\omega _1)^2 )+ 3(E-4e_1)(\wp(x) +\wp(x+\omega _1)) \\
& + (E-16e_1)(\wp(x+\omega _2) +\wp(x+\omega _3))+E^2-10e_1E-51e_1^2-\frac{9}{2}g_2, 
\end{array}$\\
$ \; \;  a(E)=8(E-7e_1), \quad c(E)=E^2-10e_1E-51e_1^2-3g_2, \quad H^{(0)} (E)=E-16e_1,$\\
$ \; \;  H^{(1)}(E)=E^4-4e_1E^3-2(7g_2+57e_1^2)E^2-4(105e_1^3+31g_3)E+81g_2^2-1114g_3e_1-127e_1^4 , $\\
$ \; \;  Q(E)= H^{(0)}(E) H^{(1)}(E), \quad  H^{(2)} (E)=1, \quad   H^{(3)} (E)=1,$\\
$ \; \;  Ht ^{(0)} (E)=(E+11e_1)(E^2-14e_1E+49e_1^2-27g_2),$\\
$ \; \;  Ht ^{(1)} (E)=E^2+22e_1E-851e_1^2+81g_2 ,$\\
$ \begin{array}{cl} Ht ^{(2)} (E)= &E^4+(12e_3+44e_2)E^3+(372e_2e_3+174e_3^2-18e_2^2)E^2 \\
& +28(e_2+e_3)(-191e_2^2-10e_2e_3+73e_3^2)E \\
& +(63132e_3e_2^3-60772e_3^3e_2+7062e_2^2e_3^2-27495e_3^4+34457e_2^4) 
\end{array}$\\
$ \begin{array}{cl} Ht ^{(3)} (E)= & E^4+(12e_2+44e_3)E^3+(372e_2e_3+174e_2^2-18e_3^2)E^2 \\
& +28(e_2+e_3)(-191e_3^2-10e_2e_3+73e_2^2)E \\
& +(63132e_2e_3^3-60772e_2^3e_3+7062e_2^2e_3^2-27495e_2^4+34457e_3^4)
\end{array}$\\
$ \; \;  H\theta (E)=E^2-\frac{62}{7}e_1E-\frac{557}{7}e_1^2-\frac{81}{7}g_2.$

\subsubsection{The case $(l_0, l_1 , l_2, l_3) =(2,2,2,2)$} $ $\\
$ \begin{array}{cl} \Xi (x,E)=& 9(\wp(x) ^2+\wp (x+\omega _1)^2 +\wp (x+\omega _2)^2 +\wp (x+\omega _3)^2) \\
& +3E(\wp(x) +\wp (x+\omega _1)+\wp (x+\omega _2) +\wp (x+\omega _3))+E^2-27g_2, 
\end{array}$\\
$ \; \;  a(E)=12E, \quad c(E)=E^2-24g_2, $\\
$ \; \;  Q(E)=(E-12e_1)(E-12e_2)(E-12e_3)(E^2-48g_2), \quad  H^{(0)}(E)=Q(E), $\\
$ \; \;  Ht ^{(0)} (E)=(E+24e_1)(E+24e_2)(E+24e_3), \quad H^{(k)} (E)=1, \quad (k=1,2,3), $\\
$ \begin{array}{cl} Ht ^{(k)} (E)= &E^6+72e_kE^5+324(-8e_k^2+g_2)E^4-3456(g_2e_k+g_3)E^3 \\
& +31104(8g_2e_k^2-4g_3e_k-g_2^2)E^2+746496(8g_3g_2e_k+4g_3^2+g_2^3-8g_2^2e_k^2), 
\end{array}$\\
$ \; \; H\theta (E)= E^6-\frac{1764}{11}g_2E^4+\frac{13824}{11}g_3E^3+\frac{72576}{11}g_2^2E^2-\frac{497664}{11}g_2g_3E-\frac{746496}{11}(g_2^3+4g_3^2).$\\

\subsubsection{The case $(l_0, l_1 , l_2, l_3) =(1,1,1,0)$}$ $\\
$ \; \; \Xi (x,E)=(E-3e_3)\wp(x)+(E-3e_2)\wp(x +\omega _1)+(E-3e_1)\wp(x +\omega _2) +E^2-\frac{3}{2}g_2, $\\
$ \; \;  a(E)=3E, \quad c(E)=E^2-\frac{3}{2}g_2 , \quad  Q(E)= \prod_{i=0}^3 H^{(i)}(E) , $\\
$ \; \; H^{(0)}(E)=E^2-3g_2, \quad H^{(1)} (E)=E-3e_1,  \quad H^{(2)} (E)=E-3e_2,  \quad H^{(3)} (E)=E-3e_3, $\\
$ \; \; Ht ^{(0)} (E)=1 , \quad Ht ^{(1)} (E)=E+6e_1 , \quad Ht ^{(2)} (E)=E+6e_2 , \quad Ht ^{(3)} (E)=E+6e_3,$\\
$ \; \; H\theta (E)=1 .$\\
Hence the functions except for $\Xi (x,E)$ are the same as the ones in the list for the case $(l_0, l_1 , l_2, l_3) =(2,0,0,0)$. Especially the monodromies of the solutions to Eq.(\ref{InoEF}) for the case $(l_0, l_1 , l_2, l_3) =(2,0,0,0)$ and the one for the case $(l_0, l_1 , l_2, l_3) =(1,1,1,0)$ are exactly the same.
For the case $(l_0, l_1 , l_2, l_3) =(0,1,1,1)$, we have \\
$\Xi (x,E)=(E-3e_1)\wp(x +\omega _1)+(E-3e_2)\wp(x +\omega _2)+(E-3e_3)\wp(x +\omega _3) +E^2-\frac{3}{2}g_2$ and the rest of the functions are the same as the ones in the list.

\subsection{The case $g=3$} \label{ssec:g3}
If the genus of the related hyperelliptic curve $\nu ^2=-Q(E)$ is three and the condition $l_0 \geq l_1 \geq l_2 \geq l_3$ is satisfied, then we have 14 cases 
$$
(l_0, l_1 , l_2, l_3) =
\left\{ 
\begin{array}{l}
(3,0,0,0), \; (3,1,0,0), \; (3,1,1,0), \; (3,1,1,1), \; (3,2,0,0), \\
(3,2,1,0), \; (3,2,2,1), \; (3,3,0,0), \; (3,3,1,1), \; (3,3,2,2), \\
(3,3,3,3), \; (2,1,1,1), \; (2,2,1,0), \; (2,2,2,0).
\end{array}
\right.
$$
We omit some twisted Heun polynomials because we do not need to use them to express the covering map.

\subsubsection{The case $(l_0, l_1 , l_2, l_3) =(3,0,0,0)$} $ $ \\
$\; \; \Xi (x,E)= 225\wp(x)^3+45E\wp(x)^2+6(E^2-\frac{75}{8}g_2)\wp(x)+E^3-15g_2 E-\frac{225}{4}g_3, $\\
$ \; \; a(E)=6(E^2-\frac{15}{4}g_2), \quad c(E)=E^3-\frac{45}{4}g_2 E-\frac{135}{4}g_3, $\\
$ \; \; H^{(0)}(E)=E, \quad H^{(k)} (E)=E^2+6e_kE+45e_k^2-15g_2, \; (k=1,2,3), $\\
$ \; \; Q(E)= \prod_{i=0}^3 H^{(i)}(E) , \quad  Ht ^{(0)} (E)=E^2-\frac{75}{4}g_2 ,$\\
$ \; \; Ht ^{(k)} (E)=E^2+15e_kE+\frac{75}{4}(g_2-12e_k^2), \; (k=1,2,3), \quad H\theta (E)=1.$

\subsubsection{The case $(l_0, l_1 , l_2, l_3) =(3,1,0,0)$} $ $ \\
$ \begin{array}{cl} \Xi (x,E)= & 225\wp (x)^3+45(E-2e_1)\wp (x)^2+(6E^2-24Ee_1-201e_1^2-\frac{75}{2}g_2)\wp (x) \\
& +(E-8e_2-3e_3)(E-3e_2-8e_3)\wp (x+\omega _1) \\
& +E^3-7e_1E^2-104e_1^2E-\frac{25}{4}g_2E-72e_1^3-\frac{225}{2}g_3. 
\end{array} $\\
$ \; \; a(E)=7E^2-13e_1E-152e_1^2-10g_2,$\\
$ \; \; c(E)= E^3-7e_1E^2-(104e_1^2+\frac{5}{2}g_2)E-102e_1^3-\frac{165}{2}g_3, $\\
$ \; \; Q(E)= \prod_{i=0}^3 H^{(i)}(E) ,\quad H^{(0)} (E)=E^2-16e_1E-32e_1^2+5g_2,$\\ 
$ \; \; H^{(1)}(E)=E^3-9e_1E^2-(117e_1^2+4g_2)E+69e_1^3-188g_3,  $\\
$ \; \;  H^{(2)} (E)=E-3e_2-8e_3, \quad H^{(3)} (E)=E-8e_2-3e_3, $\\
$ \; \;  Ht ^{(0)} (E)=E^2-\frac{58}{7}e_1E-\frac{737}{7}e_1^2-\frac{100}{7}g_2, \quad Ht ^{(1)} (E)=E^2+11e_1E-626e_1^2+50g_2, $\\
$ \; \; H\theta (E)=E-17e_1.$

\subsubsection{The case $(l_0, l_1 , l_2, l_3) =(3,1,1,0)$} $ $\\
$ \begin{array}{cl} \Xi (x,E)= & 225\wp (x) ^3+45(2e_3+E)\wp (x) ^2+\frac{3}{4}(8E^2+32Ee_3+332e_3^2-175g_2)\wp (x) \\
& + (E^2+(8e_1-10e_2)E+20e_1^2-84e_1e_2-79e_2^2 )\wp (x+\omega _1) \\
& + (E^2+(-10e_1+8e_2)E-79e_1^2-84e_1e_2+20e_2^2 )\wp (x+\omega _2) \\
& +E^3+7E^2e_3+(145e_3^2-59g_2)E+\frac{69}{2}e_3g_2-\frac{555}{2}g_3 
\end{array} $\\
$ \; \; a(E)=8E^2+26Ee_3+140e_3^2-85g_2, \quad Q(E)= \prod_{i=0}^3 H^{(i)}(E)$,\\
$ \; \; c(E)= E^3+7e_3E^2+145e_3^2E-\frac{221}{4}g_2E+42e_3g_2-255g_3$, \\
$ \; \; H^{(0)} (E)=E^3+12e_3E^2+4(60e_3^2-17g_2)E+32(2e_3g_2-5g_3),$\\
$ \; \; H^{(1)} (E)=E^2+(8e_2-10e_1)E+20e_2^2-84e_2e_1-79e_1^2, $\\
$ \; \; H^{(2)} (E)=E^2+(8e_1-10e_2)E+20e_1^2-84e_1e_2-79e_2^2,  \quad H^{(3)}(E)=1,$\\
$ \; \; Ht ^{(0)} (E)=E^2+\frac{31}{4}e_3E+\frac{323}{2}e_3^2-\frac{625}{8}g_2,$\\
$ \begin{array}{cl} Ht ^{(3)} (E)=& 
E^4+47e_3E^3+(243e_3^2-\frac{85}{4}g_2)E^2+80(7g_3-10g_2e_3)E \\
& +19332e_3^2g_2-3125g_2^2-33908g_3e_3,
\end{array} $\\
$ \; \; H\theta (E)=E^2+\frac{88}{7}e_3E+\frac{1648}{7}e_3^2-\frac{500}{7}g_2.$

\subsubsection{The case $(l_0, l_1 , l_2, l_3) =(3,1,1,1)$} $ $\\
$ \begin{array}{cl} \Xi (x,E)= & 225\wp (x)^3+45E\wp (x)^2+6(E^2-\frac{75}{4}g_2)\wp (x)+E^3-\frac{195}{4}g_2E-\frac{1125}{2}g_3  \\
& +(E-15e_2)(E-15e_3) \wp (x+\omega _1)+(E-15e_1)(E-15e_3) \wp (x+\omega _2) \\
& +(E-15e_1)(E-15e_2) \wp (x+\omega _3) ,
\end{array} $\\
$ \; \; a(E)=9(E^2-15g_2), \quad c(E)=E^3-45g_2E-540g_3, \quad Q(E)= \prod_{i=0}^3 H^{(i)}(E) ,$\\
$ \; \; H^{(0)}(E)=E^4-54g_2E^2-864g_3E-135g_2^2 , \quad Ht ^{(0)} (E)=E^2-75g_2 ,$ \\
$ \; \; H^{(k)} (E)=E-15e_k, \; (k=1,2,3),$\\
$ \; \; Ht ^{(k)} (E)=E^4+48e_kE^3+54(g_2-8e_k^2)E^2-8640e_k^3E-2025(3g_2^2+48g_3e_k-64e_k^4),$\\
$ \; \; H\theta (E)=E^3-63g_2E-540g_3.$

\subsubsection{The case $(l_0, l_1 , l_2, l_3) =(3,2,0,0)$} $ $\\
$ \begin{array}{cl} \Xi (x,E)= & 225\wp (x) ^3+45(E-6e_1)\wp (x) ^2+3(2E^2-24Ee_1-153e_1^2)\wp (x) \\
& +9(E+9e_1)\wp (x+\omega _1) ^2 + 3(E-12e_1)(E+9e_1) \wp (x+\omega _1) \\
& +E^3-21e_1E^2+9(\frac{5}{4}g_2-23e_1^2)E+9(83e_1^3-\frac{45}{4}g_3) 
\end{array} $\\
$ \; \; a(E)=9E^2-81e_1E+27(\frac{5}{4}g_2-29e_1^2),$\\
$ \; \; c(E)= E^3-21e_1E^2+9(\frac{7}{4}g_2-23e_1^2)E+9(76e_1^3-7g_3), $\\
$ \; \; Q(E)= \prod_{i=0}^3 H^{(i)}(E), \quad H^{(0)}(E)=1, \quad H^{(1)} (E)=E+9e_1, $\\
$ \begin{array}{cl} H^{(2)} (E) = & E^3+(24e_3+27e_2)E^2+(-144e_3^2-288e_2e_3-45e_2^2)E \\
& -1152e_3^3-3024e_2e_3^2-2808e_2^2e_3-711e_2^3, 
\end{array} $\\
$ \begin{array}{cl} H^{(3)} (E)= & E^3+(24e_2+27e_3)E^2+(-144e_2^2-288e_2e_3-45e_3^2)E \\
& -1152e_2^3-3024e_3e_2^2-2808e_3^2e_2-711e_3^3,
\end{array} $\\
$ \begin{array}{cl} Ht ^{(0)} (E)= & E^4-50e_1E^3+(\frac{15}{2}g_2+339e_1^2)E^2+(8208e_1^3+\frac{855}{2}g_3)E \\
& -\frac{9375}{16}g_2^2-\frac{19215}{2}g_3e_1-27153e_1^4, 
\end{array} $\\
$ \begin{array}{cl} Ht ^{(1)} (E) & =E^4-18e_1E^3+(\frac{351}{2}g_2-2133e_1^2)E^2+(38232e_1^3+\frac{9207}{2}g_3)E \\
& +\frac{50625}{16}g_2^2+\frac{60669}{2}g_3e_1-217809e_1^4,
\end{array} $\\
$ \; \; H\theta (E)=E-\frac{57}{2}e_1.$

\subsubsection{The case $(l_0, l_1 , l_2, l_3) =(3,2,1,0)$} $ $\\
$ \begin{array}{cl} \Xi (x,E)= & 225\wp (x) ^3+45(E-2e_2-6e_1)\wp (x) ^2 \\
& +3(2E^2-8(3e_1+e_2)E-128e_1^2-42e_2^2+73e_2e_1)\wp (x) \\
& +9(E+4e_1-12e_2)\wp (x+\omega _1) ^2 + 3(E+4e_1-12e_2)(E-10e_1+2e_2) \wp (x+\omega _1) \\
& +(E^2+2(e_2-15e_1)E+14e_2e_1+57e_2^2+89e_1^2)\wp (x+\omega _2) +E^3-7(e_2+3e_1)E^2 \\
&  +(-118e_1^2+212e_1e_2-15e_2^2)E+9(26e_2e_1^2-35e_2^3+89e_1e_2^2+42e_1^3) 
\end{array} $\\
$ \; \; a(E)=10E^2-4(13e_2+30e_1)E-6e_2^2-280e_1^2+752e_2e_1, $\\
$ \; \; c(E)= E^3-7(e_2+3e_1)E^2+(230e_1e_2+3e_2^2-100e_1^2)E-381e_2^3+567e_1e_2^2+300e_1^3, $\\
$ \; \; Q(E)= \prod_{i=0}^3 H^{(i)}(E), \quad H^{(0)} (E)=E+4e_1-12e_2 , \quad  H^{(1)} (E)=1,$\\
$ \begin{array}{cl} H^{(2)}(E)= & E^4-4(4e_1+e_2)E^3+(152e_1e_2-280e_1^2-10e_2^2)E^2 \\
& +(1648e_1^2e_2-740e_2^3+1792e_1e_2^2-320e_1^3)E \\
& -3279e_2^4+5520e_1^4+352e_2e_1^3+2872e_1e_2^3+8696e_1^2e_2^2,
\end{array} $\\
$ \; \; H^{(3)} (E)=E^2+2(e_2-15e_1)E+14e_2e_1+57e_2^2+89e_1^2, $\\
$ \begin{array}{cl} Ht ^{(0)} (E)= & E^4-(44e_1+\frac{12}{5}e_2)E^3+(\frac{348}{5}e_1e_2+\frac{504}{5}e_1^2-\frac{714}{5}e_2^2)E^2 \\
& +(-\frac{14016}{5}e_1^2e_2-748e_2^3+1836e_1e_2^2+\frac{51472}{5}e_1^3)E \\
& -\frac{87123}{5}e_2^4-\frac{87536}{5}e_1^4+64080e_2e_1^3-\frac{83612}{5}e_1e_2^3+27960e_1^2e_2^2,
\end{array} $\\
$ \begin{array}{cl} Ht ^{(2)} (E)= &E^3+(43e_2-37e_1)E^2+(956e_1^2-330e_1e_2-689e_2^2)E \\
& +5973e_2^3-17608e_2e_1^2+29263e_1e_2^2-18020e_1^3, 
\end{array} $\\
$ \; \; H\theta (E)=E^2-(\frac{56}{3}e_1+\frac{22}{3}e_2)E-\frac{772}{3}e_1^2+124e_2e_1-\frac{293}{3}e_2^2 .$

\subsubsection{The case $(l_0, l_1 , l_2, l_3) =(3,2,2,1)$} $ $\\
$ \begin{array}{cl} \Xi (x,E)= & 225\wp (x) ^3+45(E+4e_3)\wp (x) ^2+6(E^2+8e_3E+91e_3^2-\frac{175}{4}g_2)\wp (x) \\
& +9(E+e_1-24e_2)\wp (x+\omega _1) ^2 + 3(E+e_1-24e_2)(E-6e_1+4e_2) \wp (x+\omega _1) \\
& +9(E+e_2-24e_1)\wp (x+\omega _2) ^2 + 3(E+e_2-24e_1)(E-6e_2+4e_1) \wp (x+\omega _2) \\
& +(E+e_1-24e_2)(E-24e_1+e_2)\wp (x+\omega _3) \\
& +E^3+14e_3E^2+(349e_3^2-\frac{553}{4}g_2)E-\frac{261}{4}e_3g_2-\frac{1557}{2}g_3 
\end{array} $\\
$ \; \; a(E)=13E^2+146e_3E+2341e_3^2-760g_2, \quad Q(E)= \prod_{i=0}^3 H^{(i)}(E),$\\
$ \; \; c(E)= E^3+14e_3E^2+(349e_3^2-133g_2)E-33e_3g_2-756g_3, $\\
$ \; \; H^{(0)} (E)=1 , \quad  H^{(1)} (E)=E+e_1-24e_2, \quad   H^{(2)} (E)=E-24e_1+e_2,$\\
$ \begin{array}{cl} H^{(3)}(E)=& E^5+5e_3E^4+(-164g_2+250e_3^2)E^3+(\frac{725}{2}g_2e_3-\frac{4531}{2}g_3)E^2 \\
& +(-\frac{11227}{4}g_2e_3^2+\frac{45653}{4}g_3e_3+1216g_2^2)E \\
& +\frac{55521}{16}g_2^2e_3+\frac{165553}{4}g_3e_3^2+\frac{23777}{16}g_3g_2, 
\end{array} $\\
$ \begin{array}{cl} Ht ^{(0)} (E)= & E^6+\frac{222}{13}e_3E^5+(-\frac{5196}{13}g_2+\frac{7971}{13}e_3^2)E^4+(\frac{48953}{13}g_3-\frac{31735}{13}e_3g_2)E^3 \\
& +(\frac{3472323}{52}g_3e_3-\frac{5186397}{52}g_2e_3^2+36480g_2^2)E^2 \\
& +(-\frac{5102001}{104}g_2^2e_3-\frac{26698161}{104}g_3g_2+\frac{33717423}{26}g_3e_3^2)E \\
& +\frac{3200000}{13}g_2^3-\frac{51057923}{208}g_2^2e_3^2-\frac{1404511187}{208}g_3^2-\frac{49594475}{104}g_3g_2e_3,
\end{array} $\\
$ \begin{array}{cl} Ht ^{(3)} (E)= & E^4+100e_3E^3+(-1914e_3^2+376g_2)E^2-(8519g_3+5271g_2e_3)E \\
& +\frac{2111393}{4}g_2e_3^2-80000g_2^2-\frac{2594111}{4}g_3e_3, 
\end{array} $\\
$ \; \;   H\theta (E)=E^3+3e_3E^2+(195e_3^2-208g_2)E+\frac{2433}{4}e_3g_2+\frac{4801}{4}g_3 .$

\subsubsection{The case $(l_0, l_1 , l_2, l_3) =(3,3,0,0)$} $ $\\
$ \begin{array}{cl}  \Xi (x,E)= & 225(\wp (x)^3 + \wp (x+\omega _1)^3)+45(E-12e_1)(\wp (x)^2 + \wp (x+\omega _1)^2) \\
& +(6E^2-144e_1E+\frac{225}{4}g_2-486e_1^2)(\wp (x) + \wp (x+\omega _1)) \\
& +E^3-42e_1E^2+(\frac{75}{2}g_2-99e_1^2)E+2808e_1^3+225g_3, 
\end{array} $\\
$ \; \; a(E)= 12(E^2-24e_1E-81e_1^2+15g_2) ,\quad Q(E)=H ^{(0)} (E) H ^{(1)} (E) $\\
$ \; \; c(E)= E^3-42e_1E^2+9(5g_2-11e_1^2)E+2448e_1^3+360g_3 ,$\\
$ \; \; H^{(0)} (E)=(E-12e_1)(E^2-36e_1E+60g_2-432e_1^2) , $\\
$ \begin{array}{cl} H^{(1)}(E)= & E^4-36e_1E^3+66(g_2-9e_1^2)E^2+(4860e_1^3+324g_3)E \\
& +225g_2^2+6966g_3e_1+27297e_1^4 ,
\end{array} $\\
$ \; \; H^{(2)} (E)=1 , \quad  H^{(3)} (E)=1,$\\
$ \; \; Ht ^{(0)} (E)= (E^2-24e_1E-981e_1^2+75g_2)(E^2-54e_1E+729e_1^2-75g_2) $\\
$ \begin{array}{cl} Ht ^{(1)} (E)= & E^4-18e_1E^3+(525g_2-6291e_1^2)E^2+(173448e_1^3+21600g_3)E \\
& +22500g_2^2+153900g_3e_1-2007504e_1^4, 
\end{array} $\\
$ \begin{array}{cl} H\theta (E)= & E^4-\frac{918}{11}e_1E^3+(\frac{375}{11}g_2+\frac{16119}{11}e_1^2)E^2+(\frac{122688}{11}e_1^3+\frac{27000}{11}g_3)E \\
& -\frac{22500}{11}g_2^2-\frac{607500}{11}g_3e_1-\frac{2111184}{11}e_1^4 .
\end{array} $\\

\subsubsection{The case $(l_0, l_1 , l_2, l_3) =(3,3,1,1)$}  $ $\\
$ \begin{array}{cl}  \Xi (x,E)= & 225(\wp (x)^3 + \wp (x+\omega _1)^3)+45(E-10e_1)(\wp (x)^2 + \wp (x+\omega _1)^2)  \\
& +(6E^2-\frac{75}{4}g_2-525e_1^2-120Ee_1)(\wp (x) + \wp (x+\omega _1))\\
& + ( E^2-50e_1E+25(13e_1^2+g_2)) (\wp (x+\omega _2) +\wp (x+\omega _3)) \\
& +E^3-35e_1E^2+(\frac{5}{2}g_2-200e_1^2)E+300(13e_1^3-2g_3)
\end{array} $ \\
$ \; \; a(E)= 14E^2-340e_1E-400e_1^2+80g_2 ,$\\
$ \; \;  c(E)= E^3-35e_1E^2+10(g_2-20e_1^2)E+240(15e_1^3-2g_3)  ,\quad Q(E)=H ^{(0)} (E) H ^{(1)} (E)$\\
$ \; \;  H^{(0)}(E)=(E^2+4e_1E-80e_1^2-20g_2)(E^3-24e_1E^2+16(g_2-39e_1^2)E-1280(e_1^3+g_3)),$ \\
$ \; \;  H^{(1)} (E)=E^2-50e_1E+25(13e_1^2+g_2), \quad  H^{(2)} (E)=1 , \quad  H^{(3)} (E)=1,$\\
$ \; \; Ht ^{(0)} (E)= (E^2-50e_1E+200(5e_1^2-g_2))(E^2-\frac{80}{7}e_1E+\frac{400}{7}(g_2-23e_1^2)) $\\
$ \begin{array}{cl} Ht ^{(1)} (E)= & E^6+48e_1E^5+12(-796e_1^2+67g_2)E^4+320(33g_3+188e_1^3)E^3 \\
& +9600(-g_2^2+27g_3e_1+334e_1^4)E^2+384000(51g_3e_1^2-64e_1^5-5g_3g_2)E\\
& +640000(351g_3e_1^3-776e_1^6+30g_3^2+g_2^3) ,
\end{array} $\\
$ \begin{array}{cl} H\theta (E)= & E^6-\frac{600}{13}e_1E^5-\frac{12}{13}(700e_1^2+117g_2)E^4+\frac{320}{13}(-27g_3+1640e_1^3)E^3 \\
& +\frac{9600}{13}(-5g_2^2-201g_3e_1+886e_1^4)E^2+\frac{384000}{13}(-37g_3e_1^2-102e_1^5+5g_3g_2)E \\
& +\frac{640000}{13}(159g_3e_1^3-872e_1^6-66g_3^2+g_2^3).
\end{array} $\\

\subsubsection{The case $(l_0, l_1 , l_2, l_3) =(3,3,2,2)$} $ $ \\
$\begin{array}{cl} 
\Xi (x,E) =&  225(\wp (x) ^3 + \wp (x+\omega _1)^3)+45(E-6e_1)(\wp (x) ^2+ \wp (x+\omega _1)^2) \\
 & +(6E^2-\frac{675}{4}g_2-459e_1^2-72Ee_1)(\wp (x) + \wp (x+\omega _1)) \\
 & + 9(E-36e_1)(\wp (x+\omega _2)^2 +\wp (x+\omega _3)^2)\\
 & + 3(E+6e_1)(E-36e_1) (\wp (x+\omega _2) +\wp (x+\omega _3)) \\
 & +E^3-21e_1E^2+(-72g_2-261e_1^2)E+3249e_1^3-1917g_3 ,
\end{array} $ 
\\
$ \; \; \textstyle a(E)= 18(E^2-18e_1E-123e_1^2-15g_2) ,\quad Q(E)=H ^{(0)} (E) H ^{(1)} (E)$ \\
$ \; \; \textstyle c(E)= E^3-21e_1E^2+(-261e_1^2-63g_2)E+2853e_1^3-1773g_3$ , \\
$ \; \; H^{(0)} (E)=E-36e_1, \quad  H^{(2)} (E)=1 , \quad  H^{(3)} (E)=1,$\\
$ \begin{array}{cl} 
H^{(1)}(E)= & E^6-6e_1E^5+(-405e_1^2-189g_2)E^4+(-4212g_3-2484e_1^3)E^3 \\
 & +(8451g_2^2-55242g_3e_1+69903e_1^4)E^2 \\
 & +(-293868g_3e_1^2-275238e_1^5+345222g_3g_2)E \\
 & +1426653g_3e_1^3-7172631e_1^6+1019871g_3^2+50625g_2^3 ,
\end{array} $\\
$ \begin{array}{cl} Ht ^{(0)} (E)= & (E^4-36e_1E^3+54(37e_1^2-13g_2)E^2+108(-67e_1^3+107g_3)E \\
 & +81(625g_2^2+7382g_3e_1-13199e_1^4))(E^4+28e_1E^3-30(g_2+55e_1^2)E^2 \\
 & +468(-61e_1^3+5g_3)E-9375g_2^2-179370g_3e_1+476289e_1^4),
\end{array} $\\
$ \begin{array}{cl} Ht ^{(1)} (E)= & E^6+138e_1E^5+(-17253e_1^2+1971g_2)E^4+(-5076g_3+45036e_1^3)E^3 \\
 & +(-388125g_2^2-4881546g_3e_1+15619311e_1^4)E^2 \\
 & +(105825204g_3e_1^2-157540950e_1^5-10226250g_3g_2)E \\
 & +4936453389g_3e_1^3-5402344167e_1^6-149900625g_3^2+31640625g_2^3,
\end{array} $ \\
$ \begin{array}{cl} H\theta (E)= & E^6-\frac{390}{17}e_1E^5-(\frac{9477}{17}e_1^2+\frac{9261}{17}g_2)E^4+(\frac{213516}{17}g_3+\frac{699084}{17}e_1^3)E^3 \\
& +(\frac{421875}{17}g_2^2-\frac{10529514}{17}e_1g_3+\frac{16361919}{17}e_1^4)E^2 \\
& -(\frac{134787564}{17}g_3e_1^2+\frac{7429158}{17}e_1^5+\frac{12656250}{17}g_3g_2)E \\
& +\frac{602204301}{17}g_3e_1^3+\frac{418505049}{17}e_1^6-\frac{1287140625}{17}g_3^2+\frac{31640625}{17}g_2^3 .
\end{array} $\\

\subsubsection{The case $(l_0, l_1 , l_2, l_3) =(3,3,3,3)$} $ $\\
$ \begin{array}{cl} \Xi (x,E)= & 225(\wp (x)^3 + \wp (x+\omega _1)^3+\wp (x+\omega _2)^3 +\wp (x+\omega _3)^3) \\
& + 45E(\wp (x) ^2+ \wp (x+\omega _1)^2+\wp (x+\omega _2)^2 +\wp (x+\omega _3)^2) \\
& +(6E^2-\frac{1575}{4}g_2)(\wp (x) + \wp (x+\omega _1) + \wp (x+\omega _2) +\wp (x+\omega _3))\\
& +E^3-195g_2E-2250g_3 ,
 \end{array}$ \\
$ \; \; a(E)= 24(E^2-60g_2) , \quad c(E)= E^3-180g_2E-2160g_3,$\\
$ \: \: Q(E)=H^{(0)}(E)=E \prod _{k=1}^3 (E^2+24e_kE+720e_k^2-240g_2) ,$\\
$ \; \; H^{(1)} (E)=H^{(2)} (E)=H^{(3)} (E)=1,$\\
$ \; \; Ht ^{(0)} (E)= (E^2-300g_2)\prod _{k=1}^3 (E^2+60e_kE+300(g_2-12e_k^2)) $,\\
$ \begin{array}{l}  Ht ^{(k)} (E) =  E^{12}+288e_kE^{11}+144(41g_2-288e_k^2)E^{10}-8640(32e_k^3+17g_3)E^9 \\
 \quad +21600(-287g_2^2-3744g_3e_k+9216e_k^4)E^8+12960000(-64e_k^5+g_2g_3+16g_3e_k^2)E^7 \\
 \quad +4320000(-82944e_k^6+629g_2^3-7830g_3^2+63936g_3e_k^3)E^6 \\
 \quad +777600000(280g_3^2e_k+1024e_k^7+17g_2^2g_3-512g_3e_k^4)E^5 \\
 \quad +972000000(7056g_3^2g_2-589g_2^4-322560g_3e_k^5+62208g_3^2e_k^2+294912e_k^8)E^4 \\
 \quad +233280000000(g_2^3g_3-1920g_3e_k^6-1248g_3^2e_k^3+2560e_k^9-148g_3^3)E^3 \\
 \quad +2332800000000(-18432g_3^2e_k^4+19g_2^5+2016g_3^3e_k-36864e_k^{10}+50688g_3e_k^7-225g_2^2g_3^2)E^2 \\
 \quad +279936000000000(-120g_3^3e_k^2+96g_3^2e_k^5-g_2^4g_3+35g_3^3g_2-2048e_k^{11}+2048g_3e_k^8)E \\
 \quad +466560000000000(-54g_3^2g_2^3+729g_3^4+g_2^6) \end{array} $ \\
$ \begin{array}{l} H\theta (E)=  E^{12}-\frac{49104}{23}g_2E^{10}+\frac{3481920}{23}g_3E^9+885600g_2^2E^8-\frac{1905120000}{23}g_3g_2E^7 \\
\quad  -\frac{4320000}{23}(137g_2^3+13986g_3^2)E^6+\frac{333590400000}{23}g_2^2g_3E^5+\frac{12636000000}{23}(2160g_3^2g_2-71g_2^4)E^4 \\
\quad  +\frac{233280000000}{23}(1188g_3^3-125g_2^3g_3)E^3+\frac{114307200000000}{23}g_2^2(g_2^3-27g_3^2)E^2 \\
\quad  +\frac{2519424000000000}{23}g_2g_3(g_2^3-27g_3^2)E-\frac{466560000000000}{23}(g_2^3-27g_3^2)^2  \end{array}$ \\

\subsubsection{The cases $(l_0, l_1 , l_2, l_3) =(2,1,1,1)$, $(2,2,1,0)$, $(2,2,2,0)$}
It is shown by calculating functions $\Xi (x,E)$, $a(E)$, $c(E)$, $Q(E)$, $H ^{(i)} (E)$, $Ht ^{(i)} (E)$ $(i=0,1,2,3)$, $H\theta (E)$ explicitly that the monodromy for the case $(l_0, l_1 , l_2, l_3) =(2,1,1,1)$ (resp.  $(l_0, l_1 , l_2, l_3) =(2,2,1,0)$, $(2,2,2,0)$) is the same as the one for the case $(l_0, l_1 , l_2, l_3) =(3,0,0,0)$ (resp.  $(l_0, l_1 , l_2, l_3) =(3,1,0,0)$, $(3,1,1,1)$).

For the case $(l_0, l_1 , l_2, l_3) =(2,1,1,1)$ we have \\
$ \begin{array}{cl} \Xi (x,E) =& 9E\wp (x)^2+3E^2\wp (x)+\sum _{k=1}^3 (E^2+6e_kE+45e_k^2-15g_2)\wp (x+\omega _k) \\
& \textstyle + E^3-12g_2E-\frac{135}{4}g_3,
\end{array}$\\
and the rest of the functions are the same as the ones in the list for the case $(l_0, l_1 , l_2, l_3) =(3,0,0,0)$.

For the case $(l_0, l_1 , l_2, l_3) =(2,2,1,0)$ we have\\
$ \begin{array}{cl} 
\Xi (x,E) =& 9(E-3e_2-8e_3)\wp (x)^2+3(E+4e_2+6e_3)(E-3e_2-8e_3)\wp (x) \\
& \textstyle +9(E-8e_2-3e_3)\wp (x+\omega _1) ^2 +3(E+6e_2+4e_3)(E-8e_2-3e_3) \wp (x+\omega _1) \\
& \textstyle + (E^2-16e_1E-32e_1^2+5g_2) \wp (x+\omega _2) \\
& \textstyle + E^3-7e_1E^2-4g_2E-104e_1^2E-\frac{297}{4}g_3-135e_1^3,
\end{array}$\\
and the rest of the functions are the same as the ones in the list for the case $(l_0, l_1 , l_2, l_3) =(3,1,0,0)$.

For the case $(l_0, l_1 , l_2, l_3) =(2,2,2,0)$ we have \\
$ \begin{array}{cl} \Xi (x,E) = & 9(E-15e_3)\wp (x)^2+3(E+6e_3)(E-15e_3)\wp (x) + E^3-\frac{189}{4}g_2E-540g_3 \\
& \textstyle + 9(E-15e_2)\wp (x +\omega _1)^2+3(E+6e_2)(E-15e_2)\wp (x +\omega _1) \\
& \textstyle + 9(E-15e_1)\wp (x +\omega _2)^2+3(E+6e_1)(E-15e_1)\wp (x +\omega _2),
\end{array}$\\
and the rest of the functions are the same as the ones in the list for the case $(l_0, l_1 , l_2, l_3) =(3,1,1,1)$.

\subsection{The case $g=4$} \label{ssec:g4}
If the genus of the related hyperelliptic curve $\nu ^2=-Q(E)$ is four and the condition $l_0 \geq l_1 \geq l_2 \geq l_3$ is satisfied, then we have 24 cases 
$$
(l_0, l_1 , l_2, l_3) =
\left\{
\begin{array}{l}
(4,0,0,0), \; (4,1,0,0), \; (4,1,1,0), \; (4,1,1,1), \; (4,2,0,0), \; (4,2,1,0), \\
(4,2,1,1), \; (4,2,2,0), \; (4,2,2,2), \; (4,3,0,0), \; (4,3,1,0), \; (4,3,2,1), \\
(4,3,3,2), \; (4,4,0,0), \; (4,4,1,1), \; (4,4,2,2), \; (4,4,3,3), \; (4,4,4,4), \\
(3,2,1,1), \; (3,2,2,0), \; (3,3,1,0), \; (3,3,2,0), \; (3,3,3,1), \; (2,2,2,1).
\end{array}
\right.
$$
If $ 2 \geq l_0  \geq l_1 \geq l_2 \geq l_3$ then the genus is no more than four, and the genus equal to four case is $(l_0, l_1 , l_2, l_3) =(2,2,2,1)$. This case is related with the case $(l_0, l_1 , l_2, l_3) =(4,0,0,0)$.
Now we consider the cases $(l_0, l_1 , l_2, l_3) =(4,0,0,0)$ and $(l_0, l_1 , l_2, l_3) =(2,2,2,1)$.

\subsubsection{The case $(l_0, l_1 , l_2, l_3) =(4,0,0,0)$} $ $\\
$ \begin{array}{cl} \Xi (x,E)= & 11025\wp (x)^4+1575E\wp (x)^3+135(E^2-\frac{49}{2}g_2)\wp (x)^2 \\
& +10(E^3-\frac{371}{8}g_2E-245g_3)\wp (x)+E^4-\frac{113}{2}g_2E^2-\frac{1925}{4}g_3E+\frac{3969}{16}g_2^2, 
\end{array}$\\
$ \; \; c(E)=E^4-\frac{181}{4}g_2E^2-\frac{1295}{4}g_3E+\frac{273}{2}g_2^2, $\\
$ \; \; a(E)=10E^3-\frac{455}{2}g_2E-875g_3,  \quad  Q(E)= \prod_{i=0}^3 H^{(i)}(E),$ \\
$ \; \; H^{(0)}(E)=E^3-52g_2E-560g_3, \quad  Ht ^{(0)} (E)=E^3-\frac{343}{4}g_2E+\frac{1715}{2}g_3,  $\\
$ \; \; H^{(k)} (E)=E^2-10e_kE-35e_k^2-7g_2, \; (k=1,2,3)$\\
$ \begin{array}{cl} Ht ^{(k)} (E)= & E^4+55e_kE^3+\frac{7}{4}(77g_2-540e_k^2)E^2-490(16e_k^3+g_3)E \\
& +343(-27g_2^2-380g_3e_k+720e_k^4), \; (k=1,2,3), 
\end{array}$\\
$ \; \; H\theta (E)=E^2-\frac{196}{3}g_2.$

For the case $(l_0, l_1 , l_2, l_3) =(1,2,2,2)$ we have \\
$ \begin{array}{cl} \Xi(x,E)= & \textstyle  E^4-\frac{95}{2}g_2E^2-\frac{1295}{4}g_3E+\frac{1323}{8}g_2^2 + (E^3-52g_2E-560g_3 )\wp (x) \\
& \textstyle  +\sum _{k=1}^3 \left( 9(E^2-10e_kE-35e_k^2-7g_2)\wp (x+\omega _k )^2 \right. \\
& \left. \textstyle +3(E^3-6e_kE^2-7g_2E-75e_k^2E+28g_3-252e_k^3)\wp (x+\omega _k ) \right) , 
\end{array}$\\
and the rest of the functions are the same as the ones in the list for the case $(l_0, l_1 , l_2, l_3) =(4,0,0,0)$.
The case $(l_0, l_1 , l_2, l_3) =(2,2,2,1)$ is obtained by parallel transformation $x \rightarrow x+ \omega _3$ from the case $(l_0, l_1 , l_2, l_3) =(1,2,2,2)$. For this case, the functions except for $\Xi (x,E)$ are the same as the ones in the list.

\appendix
\section {Elliptic functions} \label{sect:append}
This appendix presents the definitions of and the formulas for the elliptic functions.

The Weierstrass $\wp$-function, the Weierstrass sigma-function and the Weierstrass zeta-function with periods $(2\omega_1, 2\omega_3)$ are defined as follows:
\begin{align}
& \wp (z)= \frac{1}{z^2}+  \sum_{(m,n)\in \Zint \times \Zint \setminus \{ (0,0)\} } \left( \frac{1}{(z-2m\omega_1 -2n\omega_3)^2}-\frac{1}{(2m\omega_1 +2n\omega_3)^2}\right),  \\
& \sigma (z)=z\prod_{(m,n)\in \Zint \times \Zint \setminus \{(0,0)\} } \left(1-\frac{z}{2m\omega_1 +2n\omega_3}\right) \nonumber \\
& \; \; \; \; \; \; \; \; \; \; \; \; \; \; \cdot \exp\left(\frac{z}{2m\omega_1 +2n\omega_3}+\frac{z^2}{2(2m\omega_1 +2n\omega_3)^2}\right), \nonumber \\
& \zeta(z)=\frac{\sigma'(z)}{\sigma (z)}. \nonumber
\end{align}
Setting $\omega_2=-\omega_1-\omega_3$ and 
\begin{align}
& e_k=\wp(\omega_k), \; \; \; \eta_k=\zeta(\omega_k) \; \; \; \; (k=1,2,3)
\end{align}
yields the relations
\begin{align}
& e_1+e_2+e_3=\eta_1+\eta_2+\eta_3=0, \; \; \; \label{eq:Leg} \\
& \eta _1 \omega _3- \eta _3 \omega _1 = \eta _3 \omega _2- \eta _2 \omega _3 = \eta _2 \omega _1- \eta _1 \omega _2 = \pi\sqrt{-1} /2, \nonumber \\
& \wp(z)=-\zeta'(z), \; \; \; (\wp'(z))^2=4(\wp(z)-e_1)(\wp(z)-e_2)(\wp(z)-e_3), \nonumber \\
& \wp (z) - \wp (\tilde{z}) = -\frac{\sigma ( z+\tilde{z})\sigma ( z-\tilde{z})}{\sigma ( z)^2\sigma (\tilde{z})^2} \nonumber
\end{align}
The addition formula for the Weierstrass zeta function is given as 
\begin{equation}
\zeta (z+\tilde{z}) -\zeta (z) -\zeta (\tilde{z}) = \frac{1}{2} \frac{\wp '(z) -\wp '(\tilde{z})}{\wp (z) -\wp (\tilde{z})}.
\label{adzeta}
\end{equation}
The periodicity of functions $\wp(z)$, $\zeta (z)$ and $\sigma (z)$ are as follows:
\begin{align}
& \wp(z+2\omega_k)=\wp(z), \; \; \; \zeta(z+2\omega_k)=\zeta(z)+2\eta_k \; \; \; \; (k=1,2,3), \label{periods} \\
& \sigma (z+2\omega _k) = - \sigma (z) \exp (2\eta _k (z + \omega _k)), \; \; \; \frac{\sigma (z+t+2\omega _k )}{\sigma (z+2\omega _k)}= \exp(2\eta _k t) \frac{\sigma (z+t)}{\sigma (z)} \nonumber
\end{align}
The constants $g_2$ and $g_3$ are defined by
\begin{equation}
g_2=-4(e_1e_2+e_2e_3+e_3e_1), \; \; \; g_3=4e_1e_2e_3.
\end{equation}
The co-sigma functions $\sigma_k(z)$ $(k=1,2,3)$ and co-$\wp$ functions $\wp_k(z)$ $(k=1,2,3)$ are defined by
\begin{align}
& \sigma_k(z)=\exp (-\eta_k z)\frac{\sigma(z+\omega_k)}{\sigma(\omega _k)}, \; \; \; \wp_k(z) = \frac{\sigma_k(z)}{\sigma(z)}, \label{eq:sigmai}
\end{align}
and satisfy
\begin{align}
& \wp_k(z) ^2 =\wp(z)-e_k, \quad \quad \quad (k,k' =1,2,3) \label{rel:sigmai} \\
& \wp _k (z+2\omega _{k'}) = \exp (2(\eta _{k'} \omega _k -\eta _{k} \omega _{k'}) ) \wp _k (z) = (-1)^{\delta _{k,k'}} \wp _k (z). \nonumber
\end{align}

{\bf Acknowledgment}
The author would like to thank Professor Etsuro Date for fruitful discussions. 
The author is partially supported by a Grant-in-Aid for Scientific Research (No. 15740108) from the Japan Society for the Promotion of Science.

\end{document}